\documentclass[reqno,draft,intlimits]{amsart}

\usepackage[active]{srcltx}

\usepackage[reqno]{amsmath}
\usepackage{amssymb}
\usepackage[all]{xy}

\usepackage{wrapfig}

\makeatletter
\long\def\@makecaption#1#2{%
  \vskip3pt
  \sbox\@tempboxa{\small#1. #2}%
    \ifdim \wd\@tempboxa >\hsize
    \small#1. #2\par
  \else
    \global \@minipagefalse
    \hb@xt@\hsize{\hfil\box\@tempboxa\hfil}%
  \fi
  \vskip0pt}
\makeatother

\newtheorem{thm}{Theorem}[section] 
  
\newtheorem{proc}{Proposition}[section]

\newtheorem{lem}{Lemma}[section]

\newtheorem{defi}{Definition}[section]
\newtheorem{rmk}{Remark}[section]
\newtheorem{ex} {Example}[section]

\def\p{\partial}

\DeclareMathOperator{\Smbl}{Smbl} \DeclareMathOperator{\id}{id} 
\DeclareMathOperator{\smbl}{smbl} 
  
\DeclareMathOperator{\Spec}{Spec}  
 \DeclareMathOperator{\im}{im}

\DeclareMathOperator{\Df}{Diff}  
\DeclareMathOperator{\Hom}{Hom}  
 \DeclareMathOperator{\Sym}{Sym}

\DeclareMathOperator{\Diff}{Diff}   
\DeclareMathOperator{\Ob}{Ob}

\newcommand{\df}{\stackrel{\mathrm{def}}{=}}

\def\cK{\mathcal K}
\def\cL{\mathcal L}
\def\cE{\mathcal E}

\def\cJ{\mathcal J}
\def\cO{\mathcal O}

\newcommand{\cA}{\mathcal A}

\newcommand{\R}{\mathbb R}
\newcommand{\Z}{\mathbb Z}
\newcommand{\C}{\mathbb C}

\newcommand{\dF}{\mathbb F}

\newcommand{\gk}{\boldsymbol{k}}

\newcommand{\gD}{\boldsymbol{\mathfrak D}}

\newcommand{\gH}{\boldsymbol{\mathfrak h}}

\newcommand{\gm}{\boldsymbol{\mathfrak m}}

\def\rw{\rightarrow}
\def\lrw{\longrightarrow}
\def\lw{\leftarrow}
\def\llw{\longleftarrow}

\DeclareFontFamily{OT1}{wncyi}{} \DeclareFontShape{OT1}{wncyi}{m}{it}{
   <5> <6> <7> <8> <9> gen * wncyi
   <10> <10.95> <12> <14.4> <17.28> <20.74> <24.88> wncyi10
  }{}
\DeclareSymbolFont{cyrletters}{OT1}{wncyi}{m}{it} 
\DeclareSymbolFontAlphabet{\cyrmath}{cyrletters} 
\DeclareMathSymbol{\rE}{\cyrmath}{cyrletters}{003} 
\DeclareMathSymbol{\rD}{\cyrmath}{cyrletters}{068} 
\DeclareMathSymbol{\rG}{\cyrmath}{cyrletters}{017} 
\DeclareMathSymbol{\rI}{\cyrmath}{cyrletters}{073} 
\DeclareMathSymbol{\rL}{\cyrmath}{cyrletters}{076} 
\DeclareMathSymbol{\rZ}{\cyrmath}{cyrletters}{090}

\def\r{\rho}

\def\be{\begin{equation}}
\def\ee{\end{equation}}
\def\bea{\begin{eqnarray}}
\def\eea{\end{eqnarray}}

\newdimen\theight
\def \refright#1{%
             \vadjust{\setbox0=\hbox{\quad\vtop{\hsize5cm\bf\noindent #1}}%
             \theight=\ht0
             \advance\theight by \dp0    \advance\theight by \lineskip
             \kern -\theight \vbox to \theight{\rightline{\rlap{\box0}}%
             \vss}%
             }}%

\begin{document}

\begin{center}
{\LARGE \bf{LOGIC OF DIFFERENTIAL  CALCULUS \\  \vspace{0.3cm}AND  THE ZOO OF 
GEOMETRIC\\ \vspace{0.3cm}STRUCTURES} }
\footnote{.}
\end{center}
\vspace{.5cm}
\begin{center}

{\Large A.M.Vinogradov}
\end{center}

\bigskip

\begin{center}
$^{1}$
Levi-Civita Institute, 83050 Santo Stefano del Sole (AV), Italia.
\end{center}
\smallskip

\vspace{1.0cm}

\noindent {\bf Abstract.}
Since the discovery of differential calculus by Newton and Leibniz and the subsequent continuous growth of its applications to physics, mechanics, geometry, etc,  it was observed that partial derivatives in the study of  various natural problems are (self-)organized in certain structures usually called geometric. Tensors, connections, jets, etc, are commonly known examples of them. This list of classical geometrical structures is sporadically and continuously widening. For instance, Lie algebroids and BV-bracket are popular recent additions into it.

Our  goal is to show that the "zoo" of all geometrical structures has a common source in the calculus of functors of differential calculus over commutative algebras, which surprisingly comes
from a due mathematical formalization of observability mechanism in classical physics. We also use this occasion for some critical remarks and discussion of some perspectives.

\vspace{5.0cm}
\pagebreak

\tableofcontents


\section{Some motivating lessons from the history} 
The dimension of the subject we shall try to embrace in this short survey forces us to start with 
some historical and philosophical observations \footnote{Philosophy is an attempt to 
explain in terms of a natural language various aspects of the surrounding us reality tacitly assuming 
that this is possible. In spite of this is, as a rule, impossible, a directing philosophy is indispensable 
in the process of  formation of that unique scientific languages in which it can be done.} 
in spite of the expected ironic reaction of the reader.

Zeno of Elea and after him Plato, Aristotle and others Greek philosophers raised the problem
of a rational/scientific/mathematical description of motion. The famous Zeno's paradoxes were
disputed for two and a half millennia \cite{Zeno1, Zeno2, Zeno3, Zeno4} and the discussions 
are going on as well even now by involving new arguments from quantum mechanics, general
relativity (\cite{Zeno5}........) and such mathematical constructions as nonstandard analysis 
(\cite{Rob}........). It is very instructive to note that as in the famous competition of Achilles and the 
tortoise  all proposed solutions of  Zeno's paradoxes were not fast enough to reach a common 
consensus before immediate counterarguments forcing them to slide away. 

The invention of differential calculus by I. Newton and G. W. Leibniz, on the other hand, made it possible
an adequate and mathematically exact description of motion, first, in mechanics and later in {\it classical physics} in general. Highest efficiency and elegancy of this approach have led working mathematicians and
physicists to a widely diffused conviction that all these philosophical discussions are something obsolete
and not very relevant. The fact that the long and tortuous path to differential calculus was paved with
various paradoxical Zeno-like logical constructions was almost forgotten.

This long history is one of others that teach that the information a human can get via his senses cannot
be adequately explained in terms of any {\it natural language}, Greek, Latin,...,English. Moreover, a 
more detailed analysis shows that the primary role of a natural language is to transmit information
(see Appendix in \cite{N}) but not to explain. Such terms as the famous {\it ``infinitesimal"} or
{\it ``wave-particle"} come to us, like relic radiation, from the periods of formation new adequate 
languages, when their terms were coined as self-contradicting hybrids in the old language by 
reflecting inadequacy.

A very general lesson to be drawn from this long story is
\begin{quote} 
{\it Even if a problem/phenomenon is clearly seen this does not automatically  imply that the scientific 
community is in possess of the adequate language for its exact mathematical description/explanation.}
\end{quote} 
By turning back to the problem of motion or, more generally, to that of ``evolution", ``continuous change",
etc, and recognizing that differential calculus is the native language for these problematics, 
at least, in the context of classical physics the next question to be posed is about the {\it ``grammar
structure"} of this language. It is clear that such concepts as ``infinitesimals", ``limits", etc,
as descriptions of our intuitive ideas in terms of a natural language can not be used for this purpose.
On the other hand, this rather natural from philosophical point of view question needs, however, to be 
put into a more concrete context allowing its scientific analysis. To this end the following observation 
is of help.

In classical physics the state of a physical system at an instant of time is completely determined 
by readings of measuring devices of a laboratory. The role of differential calculus is then
to elaborate these data in order to predict the further evolution (``motion") of the system or any
other information about it. So, it is natural to think that a due mathematical formalization of a physical laboratory should be included into the theory for its completeness.


\section{From the observation mechanism in classical physics to differential calculus}
This section is a brief summary of \cite{N} which is our starting point.
By a classical physical laboratory $\mathcal{L}$ we mean a set of all relevant measuring devices 
whose readings completely determine states of the physical system we deal with. With two devices 
$I_1,I_2\in\cL$ one may associate its sum $I_1+I_2$. This is a (virtual) device any reading of whichk   
is the sum of corresponding readings of $I_1$ and $I_2$\footnote{The modern technology allows to
easily construct this and similar devices.}. Similarly is defined the product $I_1I_2$ of $I_1$ and $I_2$.
For $\lambda\in\R, I\in\cL, \,\lambda I$ refers to the device of the same kind as $I$ but with $\lambda$-
times modified scale. We also need a ``stupid" device denoted $\mathbb{I}$ whose reading is
constantly $1$ independently of the state of the system. The role of this device is that it allows to
shift ``zero" on the scale of $I\in\cL$ by $\lambda$ by passing from $I\in\cL$ to $I+\lambda\mathbb{I}$. 
By ``constructing" all such virtual devices we obtain a commutative algebra with the unit over $\R$. 
Real, not virtual devices, measuring devices presented in a concrete laboratory are now interpreted 
to be generators of this algebra, called the \emph{algebra of observable}\footnote{Taking into 
consideration all virtual measuring devices we guarantee, beside other, the objectiveness of this 
construction, for instance, from national units of physical quantities or various suppliers of laboratory equipments}. 

Thus we mathematically formalize the concept of a classical physical laboratory as a commutative 
unitary algebra $A$ over $\R$. In these terms an observation is the assignment  to each ``measuring
device" $a\in A$ of its ``reading" $h(a)$. By definition of $A, \,h\colon A\rw\R$ is a homomorphism 
of unitary algebras.  So, the \emph{real spectrum} of $A$, denoted by 
$$
\Spec_{\R}A\df\{\mbox{all $\R$-algebra homomorphisms} \;\;h\colon A\rw\R \},
$$
is naturally interpreted as the space of all states of the physical system that we observe. 

Similarly is defined the $\gk$-spectrum of an commutative unitary algebra $A$ over a ground field $\gk$, which will be denoted $\Spec_{\gk}A$. Supplied with Zariski's topology $\Spec_{\gk}A$ becomes a topological space. A natural base of this topology consists of subsets
$$
U_a=\{h\in\Spec_{\gk}A\,\mid\,h(a)\neq 0\}\subset\Spec_{\gk}A.
$$
A homomorphism $H\colon A_1\rw A_2$ of commutative unitary $\gk$-algebras induces a map
$$
|H|\colon \Spec_{\gk}A_2\rw\Spec_{\gk}A_1, \quad |H|(h)\df h\circ H, \,\,h\in\Spec_{\gk}A_2.
$$
$|H|$ is \emph{continuous} in Zarisski's topology. If $A=C^{\infty}(M)$ with $M$ being a smooth 
manifold, then there is a natural
map $\imath_M\colon M\rw\Spec_{\R}A, \,\,M\ni z\mapsto h_z$, where $h_z(f)\df f(z)$.

The following ``spectrum theorem" shows a complete syntony of the above formalization of
the observability mechanism in classical physics with well-established facts.
\begin{thm}\label{spc}
\begin{enumerate}
\item $\imath_M$ is a homeomorphism assuming that $M$ is supplied with the standard topology.
\item Any smooth map $F\colon M\rw N$ is of the form $F=\imath_N^{-1}\circ|H|\circ\imath_M$ for
a homomorphism of unitary algebra $H\colon C^{\infty}(N)\rw C^{\infty}(M)$ and. conversely, any
homomorphism of unitary algebra $H\colon C^{\infty}(N)\rw C^{\infty}(M)$ is of the form $H=F^*$
for a smooth map $F\colon M\rw N$.
\end{enumerate}
\end{thm}

Nevertheless, having in mind that differential calculus is the native language of classical physics
the most important test of adequateness of the proposed formalization is whether 
differential calculus is somehow encoded in it. In other words, the question is whether differential 
calculus is an aspect of commutative algebra. The positive answer comes from the following definition 
and the subsequent theorem.

Let $A$ = be an unitary commutative $\gk$--algebra and $P$, $Q$ be $A$--modules. If $a\in A$
and $\nabla\colon P\rightarrow Q$ is a $\gk$--linear map, then $\delta_a(\nabla)\colon P\rightarrow Q$
is defined by $\delta_a(\nabla)(p)=\nabla(ap)-a\nabla(p), \,p\in P$. We also put 
$\delta_{a_1,\dots,a_r}\df\delta_{a_1}\circ\dots\circ \delta_{a_r}$ and observe that 
$\delta_{a_1,\dots,a_r}\df\delta_{a_1}$ is symmetric with respect to the indices $a_i$'s.
\begin{defi}\label{defdo}
$\Delta:P\rightarrow Q$ is a linear differential operator of order $\leq m$ if it is 
$\gk$--linear and $\delta_{a_0,\dots,a_m}(\Delta)=0, \,\forall a_0,a_1,\ldots, a_m\in A$.
\end{defi}
So-defined operators preserve all general elementary properties of usual differential operators.
For instance, composition of operators of orders $\leq k$ and $\leq l$, respectively, is an operator
of order $\leq k+l$, etc.
\begin{thm}
Let $A=C^{\infty}(M), \,\gk=\R$ and $\pi_i\colon E_i\rw M, \,i=1,2,$ be vector bundles over $M$.
Then the notion of a linear differential operator from $P=\Gamma(\pi_1)$ to $Q=\Gamma(\pi_2)$ 
in the sense of definition\,\ref{defdo}  coincides with the standard ones. 
\end{thm}
{\bf Note.} In what follows differential operators (DOs) will be understood in the sense of 
definition\,\ref{defdo}. Also,  we shall use ``commutative algebra" for ``commutative associative 
unitary algebra".


\subsection{Localizability of differential operators}
One of the most important properties of DOs, in view of their role in geometry and physics, is their
\emph{localizability}. More exactly this means the following.

Recall that a {\it multiplicative subset} $S\subset A$ of a commutative algebra $A$ is a subset, which 
is closed with respect to multiplication, contains $1_A$ and does not contain $0_A$. For instance, a multiplicative set $S_U=\{a\in A\,|\,h(a)\neq 0, \,\forall h\in U\}$ is naturally associated 
with an Zariski open subset $U\subset \Spec_{\gk}A$. If all elements of $S$ are products of
of elements $s_1,\dots,s_m\in S$, called {\it generators} of $S$, then $S$ is {\it finitely generated}. 

The localization of $A$ over $S$ is the 
algebra, denoted by $S^{-1}A$, formed by {\it formal fractions} $a/s, a\in A, s\in S$, 
i.e., equivalence classes of pairs $(a,s)$ with respect to the equivalence relation:
$$
(a_1,s_1)\sim(a_2,s_2) \quad\mbox{iff there is an} \quad s\in S \quad\mbox{such that} \quad s(a_1s_2-a_2s_1)=0.
$$
Addition and multiplication of formal fractions are obvious.
There is a canonical homomorphism of unitary algebras: 
$$
\iota=\iota_A:A\rw S^{-1}A, \;\iota(a)=a/1.
$$
The localization $S^{-1}P$ of an $A$--module $P$ over $S$ is defined similarly just by substituting 
$p\in P$ for $a\in A$ in the above formulae. Elements of $S^{-1}P$ are formal fractions
$p/s, \,p\in P, \,s\in S$, and $S^{-1}P$ is an $S^{-1}A$--module with respect to multiplication
 $(a/s)(p/s'=ap/ss'$. The canonical map $\iota=\iota_P:P\rw S^{-1}P$ is defined 
 by $\iota(p)=p/1$.
 
In the sequel $S^{-1}P$ will be considered as an $S^{-1}A$--module.   We shall also use the 
shortened notation $A_S, P_S$ for $S^{-1}A, S^{-1}P$, respectively, if the context allows it.
 \begin{proc}\label{LocSpec}
 $|\iota|$ imbeds  $\Spec_{\gk} A_S$ into $\Spec_{\gk} A$  and 
 $$
 \im(|\iota|)=\{h\in\Spec_{\gk} A\,|\,h(s)\neq 0, \,\forall s\in S\}=\bigcap_{s\in S}U_s.
 $$ 
If $S$ is finitely generated, then $\Spec_{\gk} A_S$ is a Zarisski open in $\Spec_{\gk} A$.
\end{proc}
 
If $U$ is a Zarisski open, then $S_U^{-1}\!A$ (resp., $S_U^{-1}P$) is called the {\it localization} of  
$A$ (resp., an $A$--module $P$) to  $U$ and will be simply denoted  by $A_U$ (resp., $P_U$).
 \begin{proc}\label{LocAlg}
 Let $\pi:E\rw M$ be a vector bundle, $A=C^{\infty}(M)$ and $P=\Gamma(\pi)$. If $U$ is an open 
 domain in $M$ identified with its image in $\Spec_{\R}A$, then
 $$
 A_U=C^{\infty}(U) \quad \mbox{and} \quad P_U=\Gamma(\pi|_U).
 $$
 \end{proc}

Let $\Delta\colon P\rw Q)$ be a DO of order $\leq k$. Its localization $\Delta_S:S^{-1}P\rw S^{-1}Q$ is defined by the formula
\begin{equation}\label{DefLoc}
\Delta_S\left(\frac{p}{s}\right)=\sum_{i=0}^k(-1)^i{k+1\choose i+1}\frac{\Delta( s^ip)}{s^{i+1}}
\end{equation}
If $S=S_U$ we will write $\Delta_U$ for $\Delta_S$.
Now we have
\begin{proc}\label{DLoc}
\begin{enumerate}
\item Formula \emph{(\ref{DefLoc})} correctly defines $\Delta_S$.
\item In the conditions of proposition\,\ref{LocAlg} $\Delta_U$ identifies with the standard
restriction of $\Delta$ to $U$.
\end{enumerate}
\end{proc}


\subsection{Geometrization of algebras and modules}\label{gm} 
Any element $a\in A$ may be ``visualized" by the associated with it $\gk$-valued function $f_a$ 
on $\Spec_{\gk}A$ (``geometrization" of $a$): 
$$
f_a(h)\df h(a), \;h\in \Spec_{\gk}A .
$$ 
These functions form a commutative algebra $A_{\Gamma}\df \{f_a | \,a\in A\}$. A natural surjective homomorphism $\gamma_A:A\rightarrow A_{\Gamma}$ is not, generally, an isomorphism. Obviously,
$$
\ker\gamma_A=\cap_{h\in \Spec_{\gk}A}\ker\,h
$$
Elements of  $\ker\gamma_A$ are in this sense {\it ``invisible"} and by this reason are called
{\it ghosts}. To stress this fact we put $\mathrm{Ghost}(A)=\ker\,\gamma_A$. A commutative
algebra $A$ without ghosts is called {\it geometric}. Obviously, 
$\Spec_{\gk}A =\Spec_{\gk}A_{\Gamma}$. Hence $A_{\Gamma}$ is geometric and $\gamma_A$
is the \emph{geometrization} homomorphism. Geometric are algebras of smooth functions on
smooth manifolds.  

In order to make ``visible"  elements of an $A$--module $P$ consider quotient modules 
$P_h=P/\ker\,h\cdot P, \,h\in \Spec_{\gk}A$. $A$--module $P_h$ may be considered as an
$A_h$-module ($A_h=A/\ker\,h)$. Since $h$ induces an isomorphism $A_h$ and $\gk$, $P_h$ 
may be considered as a $\gk$--vector space. This way the family $\{P_h\}_{h\in \Spec_{\gk}A}$ 
of $\gk$--vector spaces is associated with $P$. If $A=C^{\infty}(M)$ and $P=\Gamma(\pi)$, 
then the fiber $\pi^{-1}(x), x\in M$, is naturally identified with $P_{h_x}$.

If $p\in P$ denote by $p_h$ the coset $[p]_h=p\mod(\ker\,h\cdot P)\in P_h$. 
This way one gets a ``section" $\sigma_p:h\mapsto p_h$ of the family  $\{P_h\}$. The totality
$P_{\Gamma}$ of such sections may be viewed either as an $A$--module or as an 
$A_{\Gamma}$-module. Indeed, $\sigma_{ap}=f_a\cdot\sigma_p:h\mapsto h(a)p_h=[ap]_h$.
By definition, the \emph{support} of $P$, denoted by $\mathrm{Supp}\,P$, is 
$$
\mathrm{Supp}\,P=\mbox{Zarisski closure of}\,\{h\in \Spec_{\gk}A|P_h\neq 0\}.
$$
Denote by  $\gamma_P$ be the canonical projection $P\rw P_{\Gamma}$. Obviously,
$$
\ker\gamma_P=\cap_{h\in \Spec_{\gk}A}(\ker\,h\cdot P).
$$ 
By the same reasons as before 
elements of  $\ker\gamma_P$ are interpreted as ``invisible" or {\it ghosts}, and we put
 $$
 \mathrm{Ghost}(P)=\ker\,\gamma_P, \quad P_{\Gamma}=P/\mathrm{Ghost}(P).
 $$ 
 An $A$--module without ghosts is called {\it geometric}. The $A$--module $P_{\Gamma}$,
 called \emph{geometrization of $P$}, is, obviously, geometric.

{\bf Warning:} there are non-geometric modules over a geometric algebra and vice versa.

The following theorem by R.Swan completes the spectral theorem.  
\begin{thm}\label{geom}
Let $A=C^{\infty}(M), \,P=\Gamma({\pi})$ with $\pi$ being a vector bundle. Then 
\begin{enumerate}
\item $P=P_{\Gamma}$ and $P$ is a finitely generated projective $A$--module;
\item if $Q$ is a finitely generated projective $A$--module, then the family 
$\{Q_{h_x}\}_{x\in M}$ of $\R$-vector spaces form a vector bundle $\alpha$ over $M$
and $Q$ is canonically isomorphic to $\Gamma(\alpha)$.
\end{enumerate}
\end{thm}

If $H:A\rw B$ is a homomorphism of commutative $\gk$--algebras, then any $B$--module $Q$ 
acquires an $A$--module structure with respect to multiplication 
$$
(a,q)\mapsto a*q\df H(a)q, a\in A, q\in Q.
$$ 
This structure will be called {\it $H$--induced}.
\begin{proc}\label{ChageGeom}
If $Q$ is a geometric $B$ module, then it is a geometric $A$-module with respect to the
$H$--induced module structure. 
\end{proc}

Finally, we emphasize that mportance of \emph{geometricity}, in particular, comes from the fact 
that the standard differential geometry is a part of differential calculus over commutative algebras 
in the category of geometric modules over smooth function algebras. 


\subsection{Basic notation}\label{BN}
Here we fix basic notation that will be used in the sequel.
Let $A$ be a commutative $\gk$-algebra and $P,Q$ be $A$-modules.
Denote the totality of all DOs $\Delta:P\rightarrow Q$ of order $\leq k$
by $\Df_k(P,Q)$. $\Df_k(P,Q)$ has a natural $A$-bimodule structure. The \emph{left}
(resp., \emph{right}) $A$-module  structure of it is defined by 
\begin{equation}\label{lr}
(a,\Delta)\mapsto a_Q\circ \Delta 
\quad (\mbox{resp.},  (a,\Delta)\mapsto \Delta\circ a_P)
\end{equation}
where $a_R$ stands for the multiplication by $a\in A$ operator in an $A$-module $R$.
$\Df_k^{<}(P,Q)$ and $\Df_k^{>}(P,Q)$ denote the corresponding left and right
$A$-modules, respectively. The totality $\Df(P,Q)$ of all DOs from $P$ to $Q$ is a filtered $A$-bimodule
$$
\Hom_A(P,Q)=\Df_0(P,Q)\subset\dots\subset\Df_k(P,Q)\subset\Df_{k+1}(P,Q)\subset\dots\Df(P,Q)
$$
$\Df(P,P)$ is a filtered $\gk$-algebra with respect to composition of DOs,
and $\Df(P,Q)$ is a left (resp., right) filtered $\Df(Q,Q)$-module (resp., $\Df(P,P)$-module).
Also, we use  the  short notation $\Df P$ for $\Df(A,P)$. Similar meaning has $\Df_k^<P$, etc.

\section{Back to Zeno and tangent vectors}\label{ZT}
The key notion of classical differential calculus, namely, that of derivative, opened the way to constructively manipulate with previously intuitive ideas of velocity, acceleration, etc, in mechanics and of tangency,
curvature, etc, in geometry. It is remarkable that the mechanism of observability allows to transform 
the ``antique greek intuition" into a rigorous definition in a very simple and direct way. Indeed, if  $h_t$ 
is the state of the mechanical system at the instant of time $t$, then its motion is described by the curve
 $\gamma\colon t\mapsto h_t$ on $\Spec_{\R}A$ with $A$ being the algebra of observables.  Intuitively,  
 the velocity of this motion at the instant $t$ is the tangent vector to $\gamma$ and, therefore, to 
 $\Spec_{\R}A$ at the point $h_t$. In the old-fashioned terms this vector should be
 $$
\xi=\lim_{\Delta\,t\rw 0}\frac{1}{\Delta\,t}(h_{t+\Delta\,t}-h_t)
$$ 
This expression, as an $\R$-linear map from $A$ to $\R$,  is well-defined assuming that the limit exists. 
It is easily deduced that the relation $\xi$ with the multiplicative structure of $A$ is described by the 
``Leibniz rule"
$$
\xi(ab)=\xi(a)b+a\xi(b), \;\forall a,b\in A.
$$ 
In the general algebraic context these heuristic considerations motivate to adopt the following definition.
\begin{defi}\label{tv}
Let $A$ be a commutative algebra over $\gk$. A linear map $\xi\colon A\rw\gk$ is a tangent vector to
$\Spec_{\gk}A$ at its point $h$ if $\xi(ab)=\xi(a)h(b)+h(a)\xi(b), \;\forall a,b\in A$.
\end{defi}
In the case $A=C^{\infty}(M), \,\gk=\R$ this definition gives standard tangent vectors to $M$ if $M$ is
identified with $\Spec_{\R}A$ via the spectrum theorem. Similarly, derivations of a commutative
$\gk$-algebra $A$ are naturally interpreted as vector fields on $\Spec_{\gk}A$, etc.

Definition \ref{tv} has the following useful interpretation. With a point $h\in\Spec_{\gk}A$ one can
associate an $A$-module structure in the field $\gk$ by defining the module product to be
$a\cdot\lambda\df h(a)\lambda, \,a\in A, \,\lambda\in \gk$. Denote this $A$-module by $\gk_h$.
Then a tangent vector to $\Spec_{\gk}A$ at the point $h$ may be viewed as a first order DO 
$\xi\colon A\rw\gk_h$ such that $\xi(1_A)=0$. This illustrates universality of definition\,\ref{defdo}.
and the way to transform spectra of commutative algebras into objects of new pithy differential geometry. 

\begin{ex} Let $N$ be a closed subset of a smooth manifold $M$. Define the smooth function algebra 
on $N$ by putting $C^{\infty}(N)\df C^{\infty}(M)|_N$. The subset $N$ supplied with the algebra will
be called a smooth set. The spectrum of the algebra $C^{\infty}(N)$ is naturally identified with $N$
and, therefore, differential geometry of $N$ can be developed along the lines as above. It is 
nontrivial even for rather exotic smooth subsets. For instance, tangent spaces to the Cantor set
are 1-dimensional and they are 2-dimensional for the Peano curve. Also, vector fields on the Cantor 
set are nontrival and trivial on the Peano curve. The algebra of differential forms (see below) is
nontrivial on these smooth sets.
\end{ex}
\begin{ex} 
All DOs of order greater then zero over the algebra of continuous functions on a smooth manifold 
$M$ are trivial, while the $\R$-spectrum of this algebra is naturally identified with $M$.
\end{ex}

The above general algebraic formalization of the intuitive idea of velocity makes it possible to ``prove" impossibility to adequately describe the phenomenon of motion in a natural language.  To this end
one has to analyze a system of statements (propositions) pretending to such a description, which  
formally takes part of a \emph{propositional} or Boolean algebra. Recall that the basic operations 
of a Boolean algebra are
conjunction $(\land)$, disjunction $(\vee)$ and negation  $(\neg)$. However, for computations of the
truth value of propositions the operations of addition $(x\oplus y\df (x\land\neg y)\vee (\neg x\land y)$
and multiplication ($x\cdot y\df x\wedge y$) are more convenient. This leads to the equivalent notion of a
Boolean ring, i.e., a commutative unitary algebra $(A, \oplus, \cdot)$ over the field $\dF_2$ of integers
modulo 2 with the property $a\cdot a=a, \,\forall a\in A$ (\emph{idempotence}). In this setting  truth 
values are interpreted to be algebra homomorphisms $A\rw \dF_2$, i.e., elements of the spectrum 
$\Spec_{\dF_2}A$. So, informally, a Boolean ring may be viewed as an $\dF_2$-algebra of observables, 
i.e., a ``laboratory" whose ``measuring devices" are supplies with the \{false, true\}-scale. A simple computation based on idempotence property of $A$ proves
\begin{proc}\label{doboole}
All DOs of order greater than zero over a Boolean ring $A$ are trivial. In particular, all tangent vectors 
to $\Spec_{\dF_2}A$ are trivial.
\end{proc}
Therefore, by adopting the observability principle we see that the phenomenon of motion is 
inexpressible in terms of a natural language. In particular, Zeno's paradoxes are not, in fact, 
paradoxical from this point of view, since the truth value of any reasoning depends on truth 
values assigned to single propositions forming this reasoning. This, however, is too personal
as one can see from texts dedicated to Zeno's paradoxes and hence can not be objectively 
resolved in terms of this language.

It should be stressed that triviality of DOs over a Boolean ring is not due to discreteness
of the ground field $\dF_2$. For instance, the spectrum of the algebra $\dF_2[x]$ of polynomials with 
coefficients in $\dF_2$ consists of two points $h_0\colon p(x)\mapsto p(0)$ and 
$h_1\colon p(x)\mapsto p(1), \,p(x)\in \dF_2[x]$. Then it is directly follows from the definition 
that there is exactly one nontrivial tangent vector $\xi_{\epsilon}$ at $h_{\epsilon}, \epsilon = 0,1$, namely, 
$\xi\colon p(x)\mapsto p'(\epsilon)$ where $p'(x)$ stands for the formal derivative of $p(x)$. 
Moreover, the vector field $X\colon h_\epsilon\mapsto\xi_{\epsilon}$ on $\Spec_{\dF_2}(\dF_2[x])$
generates the flow $A_t, \,t\in\dF_2,$ defined by 
$$
A_t^*\df e^{tX}\colon (\dF_2[x])_{\Gamma}\rw (\dF_2[x])_{\Gamma}, \,t\in\dF_2,
$$
which is well-defined since $X^2=0$. It is easy to see that $A_0=\id$ and $A_1$ interchanges
$h_0$ and $h_1$.  (Note that $\dF_2[x]$ is not geometric.)
\begin{rmk}
This simple example gives a counterexample to our intuition for which
``discreteness" and ``differential calculus" are incompatible matters.
\end{rmk}


\section{The shortest way from observability to Hamiltonian mechanics}
Here we shall show how the mathematical framework for Hamiltonian mechanics can be rediscovered
by answering a natural from the ``observability philosophy" question:
\begin{quote}
{\it What is the algebra of observables for $T^*M$ assumed that  $C^{\infty}(M)$ is the algebra
of observables for $M$?}
\end{quote} 

First, we associate with filtered modules and algebras $\Df(P,Q), \Df\,P$, etc, the corresponding graded objects called (main) \emph{symbols}:
$$                   
\Smbl_k(P,Q)=\frac{\Df_k(P,Q)}{\Df_{k-1}(P,Q)}, \quad k\geq 0,
$$
assuming that $\Df_{-1}(P,Q)=0$, and
$$
\Smbl(P,Q)=\bigoplus_{k\geq 0}\Smbl_k(P,Q). 
$$

If $\Delta\in\Df_k(P,Q)$,
then $\smbl_k\Delta\df(\Delta \;\mathrm{mod}\;\Df_{k-1})\in\Smbl_k(P,Q)$ is called the (main)
symbol of $\Delta$. The composition of DOs induces the {\it composition
of symbols}
\begin{equation}\label{SMBLcomposition}
\Smbl_k(P,Q)\times\Smbl_l(Q,R)\stackrel{composition}{\longrightarrow}\Smbl_{k+l}(P,R).
\end{equation}

The induced from $\Df_k(P,Q)$ left and right $A$-module structures on  $\Smbl_k(P,Q)$
coincide, since $\delta_a(\Df_k(P,Q))\subset\Df_{k-1}(P,Q)$. So, the composition of
symbols (\ref{SMBLcomposition}) is $A$-bilinear. In particular, $\Smbl(R,R)$ is an 
associative graded
$A$-algebra, and $\Smbl(P,Q)$ is a right graded $\Smbl(P,P)$-module and a left graded
$\Smbl(Q,Q)$-module. As in the case of differential operators we shall use $\Smbl P$ for
$\Smbl(A,P)$. In particular, $\Smbl A$ is a graded $A$--algebra and $\Smbl\,P$ is a
graded $\Smbl\,A$--module.

The natural Hamiltonian formalism is based on the following elementary lemma.
\begin{lem}\label{SymCom}
If $A$--module $P$ is 1-dimensional, then the algebra $\Smbl(P,P)$ is commutative and 
hence $[\Delta,\nabla]\in\Df_{k+l-1}(P,P)$ if $\Delta\in\Df_k(P,P)$ and $\nabla\in\Df_l(P,P)$.
\end{lem}

This allows to define the {\it Poisson bracket} in $\Smbl(P,P)$ for
an 1-dimensional $A$--module $P$ by putting:
\begin{equation}\label{PoiBra}
\{\smbl_k\Delta,\smbl_l\nabla\}\df\smbl_{k+l-1}[\Delta,\nabla]
\end{equation}
This bracket inherits skew-commutativity and the Jacobi identity from the commutator 
$[\cdot,\cdot]$. So, $(\Smbl(P,P),\{\cdot,\cdot\})$ is a Lie algebra over $\gk$. Moreover, 
$\{\cdot,\cdot\})$ is a \emph{bi-derivaation} as it follows from the elementary property
$[\Delta\circ\nabla, \square]=\Delta\circ[\nabla,\square]+[\Delta,\square]\circ\nabla$ of
 commutators. So, $X_s\df\{s,\cdot\}, \,s\in\Smbl(P,P)$, is a derivation of $\Smbl(P,P)$
 and hence may be viewed as a vector field on $\Spec_{\gk}(\Smbl(P,P))$, which will be 
called \emph{Hamiltonian}.
\begin{thm}\label{ham}
\begin{enumerate}
\item The algebra $\smbl\,C^{\infty}(M)$ is geometric.
\item $\Spec_{\R}(\Smbl\,C^{\infty}(M))$ is canonically identified with the total space of the
cotangent bundle 
$T^*(M)$ of $M$ so that the elements of $\smbl\,C^{\infty}(M)_{\Gamma}$ are identified
with smooth functions on $T^*(M)$, which are polynomial along fibers of $T^*(M)\rw M$.
\item  The Poisson bracket in $\Smbl\,C^{\infty}(M))$ is identified with the standard Poisson 
bracket on $T^*(M)$.
\item If $f\in C^{\infty}(M)$,  then for the section $\sigma_{df}\colon M\rw T^*(M), \,x\mapsto d_xf$  
of the cotangent bundle we have
$$
\sigma_{df}^*(\smbl_k\Delta)=\frac{1}{k!}\delta_f^k(\Delta)
$$
where $\smbl_k\Delta$ is interpreted as a function on $T^*(M)$.
\end{enumerate}
\end{thm}
See \cite{N} for a proof. 

Construction of symbols is naturally localizable. Indeed, there is a natural homomorphism
$$
\iota_{smbl}:\Smbl_k(P,Q)_S\rw\Smbl_k(P_S,Q_S), \;\iota(\smbl_k\Delta)\mapsto\smbl_k\Delta_S,
\;\Delta\in\Df_k(P,Q),
$$
of $A$-modules, which respects the product of symbols and 
hence all involved module structures and, as a consequence, the Poisson bracket. 

Theorem \ref{ham} reveals the true nature of the standard Poisson and, therefore, of symplectic 
structure on $T^*(M)$ and as such has important consequences. In particular, it directly leads to 
a conceptual definition, namely,  as $\Spec_{\gk}(\Smbl A)$, of the cotangent bundle to the spectrum 
$\Spec_{\gk}A$ of a commutative algebra $A$. Its mechanical interpretation is $\colon$ if 
$\Spec_{\gk}A$ is the configuration space of a ``mechanical system", then $\Spec_{\gk}(\Smbl A)$ 
is its ``phase space". Probably, the most important consequence of this interpretation is that it
offers a direct algorithm for developing Hamiltonian mechanics/formalism over graded commutative 
algebras, for instance, over super-manifolds. This is a powerful instrument of constructing adequate mathematical models in physics, whose potential is far from being used  at large. 

The somehow mysterious role of symbols in this construction is, in fact, absolutely natural in view
of the fact that  \emph{propagation of fold-type singularities} of solutions of (nonlinear ) PDE $\cE$ 
is described by the Hamiltonian vector field with the main symbol of $\cE$ as its Hamiltonian. The corresponding Hamilton-Jacobi equation (= the ``eikonal equation" in geometrical optics) is only
one in the system that completely describes the behavior of the fold-type singularities 
\cite{L-C, Lun, Mas, Tay, Vgs,  L, eik}. Unfortunately, this fact seems not to be sufficiently known. 
Additionally, it sheds a new light on the quantization  problem \cite{Vgs, VWsym}.
\begin{ex}
Fold--type singularities of solutions of $u_{xx}-\frac{1}{c^2}u_{tt}-mu^2=0$ are described by the
following equations assuming that the wave fronts is in the form $x=\varphi(t)$ and
$y\df u|_{\textrm{wavefront}},\quad h\df u_x|_{\textrm{wave front}}$:
$$
\left\{\begin{array}{c}\ddot{y}+(cm)^2g=\pm 2c\dot{h} \\1-\frac{1}{c^2}\dot{\varphi}^2=0\Leftrightarrow\dot{\varphi}=\pm c\end{array}\right.\Leftarrow\left[\begin{array}{c}\textrm{Equations describing} \\\textrm{behaviour of fold--type} \\\textrm{singularities}\end{array}\right.
$$

In particular, if $\dot{h}=0$ (``resting particle"), then $\ddot{g}+(cm)^2g=0$ $\Rightarrow$ 
$\boxed{\nu=mc}$
\end{ex}

The following toy example illustrates another aspect of the above approach, namely,  the possibility
to develop ``discrete" Hamiltonian formalism.
\begin{ex}
Let $A=\dF_2[x]$. Then $\Spec_{\dF_2}(\Smbl A)$ consists of 4 points and Hamiltonian
vector fields on it form a 3-dimensional Lie algebra $\{e_1,e_2,e_3\}$ over $\dF_2$ with
$[e_1,e_2]=0, [e_1,e_3]=e_1, [e_2,e_3]=e_2$ \emph{(}see the end of section \emph{\ref{ZT})}.
\end{ex}

Finally, we note that by adding to the above mathematical scheme the physical principles prescribing
how to associate the energy function with a given mechanical system one gets a complete physical 
theory of mechanics in the Hamiltonian form. It is worth stressing that in this approach such 
fundamental for the traditional mechanics concepts as that of force, etc, are some consequences
of the above general principles.  Unfortunately, the textbook based on this approach is still waiting 
to be written. Finally,  by remembering the long and not very 
straight-line history of Hamiltonian mechanics one may note now that the conceptual
algebraic approach to differential calculus is a good shortcut.


\section{Basic functors of differential calculus}
In thus section we give some basic examples of \emph{functors of differential calculus} (FDC)
and of relating them natural transformations. These two ingredients form what can be informally 
called the ``logic of differential calculus". 

From now on $A$ stands for a commutative $\gk$-algebra and $P,Q,...$ for $A$-modules.


\subsection{A play between the left and the right}
The most obvious FDCs are $\Df_k^<\colon P\mapsto \Df_k^<P, \,k\geq0,$ and similarly for $\Df_k^>$ 
and $\Df_k$. They will be used to construct further FDCs and some relating them natural transformations.

We start by noticing that the operators 
$$
\id_k^<\colon \Df_k^<(P,Q)\rw \Df_k^>(P,Q), \quad\id_k^>\colon \Df_k^>(P,Q)\rw \Df_k^<(P,Q),
$$
which  are identity maps as maps of the \emph{set} $\Df_k(P,Q)$, are, generally, $k$-th order DOs. 
Indeed, it directly follows from $\delta_a(\id_k^<)=(a_Q-a_P)|_{\Df_k(P,Q)}$ and similarly for $\id_k^>$
(see (\ref{lr})).

The operator $\rD_k^<\colon  \Df_k^<P\rw P,  \,\Delta\mapsto \Delta(1_A)$ is, obviously, a 
homomorphism of $A$-modules. i.e. a 0-th order DO, while the operator
$\rD_k^>=\rD_k^<\circ \id_k^<$, which coincides with $\rD_k^<$ as the map of sets, is, generally, 
of order $k$. The operator $\rD_k^>$ is \emph{co-universal} in the following sense.

Let $\Delta\in \Df_k(P,Q)$. The homomorphism 
$$
h^{\Delta}\colon P\ni p\mapsto \Delta_p\in \Df_k^>Q, \quad \Delta_p(a)\df \Delta(ap)
$$ 
makes the diagram
\begin{equation}\label{RusD}
\xymatrix{
P \ar[r]^-{h^{\Delta}} \ar[rd]^{\Delta} & \Df_k^>Q \ar[d]^{\rD_k^>} \\
& Q
}
\end{equation}
commutative. Then it is easy to see that the map 
\begin{equation}\label{co-jet}
\Df_k^>(P,Q)\rw\Hom_A(P, \Df_k^>Q), \quad\Delta\mapsto h^{\Delta}
\end{equation}
is an isomorphism of $A$-modules, which \emph{co-represents} the functor $P\mapsto \Df_k(P,Q)$
($Q$ is fixed) in the category of $A$-modules.  Mapping (\ref{co-jet}) considered in the left $A$-module
structures is an isomorphism of $A$-modules as well:
\begin{equation}\label{co-jet1}
\Df_k^<(P,Q)\rw\Hom_A^<(P, \Df_k^>Q).
\end{equation}
The upper index $``<"$ in $\Hom_A^<(P, \Df_k^>Q)$ tells that the set of 
$A$-homomorphisms from $P$ to $\Df_k^>Q$ is supplied with the $A$-module structure
induced by the left $A$-module structure in $\Df_kQ$. By specifying (\ref{RusD}) for 
$\Delta=\id_k^<$ we have:
\begin{equation}\label{}
\Df_k^<(P,Q)\stackrel{h^{\id_k^<}}{\lrw} \Df_k^>(\Df_k^>(P,Q)) \;\mbox{and if} \;P=A, \;\Df_k^<Q\rw \Df_k^>(\Df_k^>Q)
\end{equation}

The $l$-th \emph{prolongation} $h^{\Delta}_l$ of $h^{\Delta}, \,l\geq0,$ is defined to be 
$h^{\Delta}_l\df h^{h^{\Delta}\circ \rD_l^>}$. It makes the following diagram  commutative:
\begin{equation}\label{h-pro}
\xymatrix{
\Df_l^>P \ar[r]^-{\rD_l^>} \ar[d]_{h^{\Delta}_l}& P \ar[d]^{h^{\Delta}} \ar[rd]^{\Delta} & \\
\Df_l^>(\Df_k^>Q) \ar[r]^-{\rD_l^>} & \Df_k^>Q  \ar[r]^-{\rD_k^>\;} & Q
}
\end{equation}

Since the order of DO $\Delta_{(l)}=\Delta\circ\rD_l^>$ is $\leq k+l$, the diagram
\begin{equation}\label{d-l}
\xymatrix{
 \Df_l^>P \ar[d]_{\rD_l^>} \ar[r]^-{h^{\Delta_{(l)}}} \ar[rd]^-{\Delta_{(l)}} &  \Df_{k+l}^>P \ar[d]^{\rD_{k+l}^>} \\
   P \ar[r]^-{\Delta} & Q
}
\end{equation}
whose upper triangle is (\ref{RusD}) for $\Delta_{(l)}$, is commutative. The diagram 

\begin{equation}\label{comp}
\xymatrix{
 \Df_l^>(\Df_k^>Q)  \ar[r]^-{c_{l,k}} \ar[d]_{\rD_l^>}& \Df_{k+l}^>Q \ar[d]^{\rD_{k+l}^>}\\
 \Df_k^>Q \ar[r]^-{\rD_k^>} & Q
}
\end{equation}
is the particular case of (\ref{d-l}) for $\Delta=\rD_k^>$. By definition. $c_{l,k}=h^{\square}$ with
$\square=\rD_k^>\circ\rD_l^>\colon  \Df_l^>(\Df_k^>Q)\rw Q $.

Now, by combining diagrams (\ref{h-pro})-(\ref{comp}), we get the commutative diagram
\begin{equation}\label{h}
\xymatrix{
 \Df_l^>P \ar[r]^-{\rD_l^>} \ar[d]_{h^{\Delta_{(l)}}} \ar[dr]^-{\Delta_{(l)}}& P \ar[d]^{\Delta} \\
\Df_{k+l}^>Q  \ar[r]^-{\rD_{k+l}^>}\; & Q
}
\end{equation}
and the homomorphism of filtered $A$-modules 
$$
h^{\Delta}_*\colon \Df^>P\rw\Df^>Q, \quad h^{\Delta}_*|_{ \Df_l^>P}=h^{\Delta_{(l)}}.
$$
Obviously, $h^{\Delta}_*(\square)=\Delta\circ\square,\,\square\in\Df^>P$.

Representing functors $\Df_k$ objects are \emph{jets} (see section\,\ref{rep}) and some 
constructions with them come from diagrams (\ref{h-pro})-(\ref{h}). 
\begin{ex}\label{right}
The symbol $^<\Df_1^> A\otimes_A P$ refers to the $A$-module that as a $\gk$-vector space 
coincides with $\Df_1^> A\otimes_A P$, while its $A$-module structure is induced by the left
multiplication in $\Df_1A$.
One of possible definitions of a connection in an $A$-module $P$ is an $A$-homomorphism
$\kappa\colon ^<\Df_1^> A\otimes_A P\rw P$ such that $\kappa(\id_A\otimes p)=p$. In this 
terms the covariant derivative $\nabla_X\colon P\rw P, \,X\in D(A)$ is defined by 
$\nabla_X(p)\df\kappa(X\otimes p)$. The ``right" analogue of it, a \emph{right connection}
in $P$, is defined as an $A$-homomorphism $\xi\colon ^>\Df_1^< A\otimes_A P\rw P$. The
\emph{right} covariant derivative $_X\nabla$ corresponding to $X\in D(A)$ is defined by
$_X\nabla(p)\df-\kappa(X\otimes p)$. Such a connection is \emph{flat} if 
$[_X\nabla,_Y\nabla]=_{[X,Y]}\nabla$. Right connections are used in the construction of integral 
forms $($see subsection\,\ref{int} and \cite{Man, MMV}$)$.
\end{ex}


\subsection{Multi-derivation functors}\label{m-der}
Now we shall introduce not so obvious functors, which correspond to multi-vector fields on manifolds
for smooth function algebras. We start from the \emph{derivation} functor $D$:  
$$
D\colon P\mapsto D(P)=\{\Delta\in\Diff_1^<(P)\ |\ \Delta(1)=0\}\equiv
$$
$$\equiv\{\Delta:A\longrightarrow P\ |\ \underset{\textrm{derivations}}{\underbrace{\Delta(ab)=
a\Delta(b)+b\Delta(a)}}\}$$
If the $A$-module $D(A)$ is geometric, then its elements are vector fields on $\Spec A$. The functor 
$D$ is followed by functors
\begin{equation}\label{D2}
D_2\colon P\mapsto D_2(P)\df D(D(P)\subset\Diff_1^>P)
\end{equation}
\begin{equation}\label{P2}
\gD_2\colon P\mapsto\gD_2^{(>)}(P)\df\Diff_1^{(>)}(D(P)\subset\Diff_1^>P)
\end{equation}

Here the symbol $D(D(P)\subset\Diff_1^>P)$ stands for the totality of all derivations 
$\Delta\colon A\rw\Diff_1^>P$ such that $\Delta(A)\subset D(A)$ supplied with the $A$-module
structure $(a, \Delta)\mapsto a_P\circ\Delta$. The symbol $\Diff_1(D(P)\subset\Diff_1^>P)$
has the similar meaning, while $\Diff_1^>(D(P)\subset\Diff_1^>P)$ refers to the same 
$\gk$-vector space but supplied with the $A$-module structure $(a, \Delta)\mapsto \Delta\circ a_A$.

Now it is important to observe that the imbedding $D_2(P)\subset \gD_2^>(P)$ of $\gk$-vector spaces
is a 1-st order DO over $A$. This makes meaningful the following inductive definition :
\begin{equation}\label{Dm}
D_m(P)\df D(D_{m-1}(P)\subset \gD_{m-1}^>(P))
\end{equation}
\begin{equation}\label{Pm}
\gD_m^{(>)}(P)=\Diff_1^{(>)}(D_{m-1}(P)\subset \gD_{m-1}^>(P))
\end{equation}
If $A=C^\infty(M)$, then $D_m(A)$ consists of $m$-vector fields on $M$. 

It is not difficult to deduce from the above definitions 
natural embeddings of $A$-modules
\begin{equation}\label{}
D_{m+n}(P)\subset D_m(D_n(P)), \quad D_m(P)\subset D_{m-1}^<(\Df_1^>(P)
\end{equation}
and the splitting $\gD_m(P)=D_{m-1}(P)\oplus D_m(P)$. In other words, there are natural 
transformations of FDOs
\begin{equation}\label{trans}
D_{m+n}\rw D_m\circ D_n,  \quad D_m\rw D_{m-1}^<(\Df_1^>), \quad D_{m-1}\lw\gD_m\rw D_m.
\end{equation}
If $A=C^{\infty}(M)$, then the module $\Lambda^m(M)$ of $m$-th order differential forms on $M$
is the representing object for the FDC $D_m)$ in the category of geometrical $A$-modules. Below 
we shall see that transformations (\ref{trans}) are ``responsible" for well-known properties of differential 
forms


\subsection{The meaning of multi-derivation functors and $\Df$-Spencer complexes}\label{Spen}
At the first glance enigmatic definitions (\ref{D2})-(\ref{Pm}), in fact, appear rather naturally as the
following diagram shows:

\begin{equation}\label{sp2}
\xymatrix{
0 \ar[r] &  D_2(P) \ar[r] \ar[d] & D^<(\Df_1^>P) \ar[r]^-{c_{1,1}} \ar[d] & 
\Df_2(P) \ar[r]^-{\rD_2} \ar[d]^= & P \ar[r] \ar[d]^{=} \ar[d] & 0\\
0 \ar[r] & \gD_2(P) \ar[r] & \Df_1^<(\Df_1^>(P)) \ar[r]^-{c_{1,1}}  & 
\Df_2(P)\ar[r]^-{\rD_2}  & P \ar[r]   & 0
}
\end{equation}
where all non labeled maps are natural embeddings. As it directly follows from the definitions
$\gD_2(P)$ is exactly the description of the kernel of $c_{1,1}$ and $D_2(P)$ of the kernel of
its restriction to $D^<(\Df_1^>P)$. Moreover, both lines in (\ref{sp2}) are complexes, and the upper
line is called the \emph{2-nd order $\Df$-Spencer complex}. Similarly, a more complicated diagram 
chasing shows that $D_m(P)$ is the kernel of a natural map 
$D_{m-1}^<(\Df_1^>P)\rw D_{m-2}^<(\Df_2^>P)$. The imbedding 
\begin{equation}\label{sp-m}
D_m(P)\hookrightarrow D_{m-1}^<(\Df_1^>P)
\end{equation} 
may be called the \emph{$(m-1)$-th exterior co-differential}, since it is represented by the exterior 
differential $d\colon \Lambda^{m-1}(M)\rw\Lambda^m(M)$ if $A=C^{\infty}(M)$ (see below). Also, 
by applying (\ref{sp-m}) to $\Df_k^>(P)$ instead of to $P$ we obtain
$$
D_m^<(\Df_k^>(P))\rw D_{m-1}^<(\Df_1^>(\Df_k^>(P))\rw D_{m-1}^<(\Df_{k+1}^>(P))
$$
where the right map is induced by $c_{1,k}\colon \Df_1^>(\Df_k^>(P)\rw \Df_{k+1}^>(P)$. This
composition is the differential in the $(m+k)$-th \emph{$\Df$-Spencer complex} 
\,$\mathrm{Sp_{m+k}}(P)$:
\begin{equation}\label{d-sp}
\xymatrix@1{
0 \ar[r] & D_{m+k}(P) \ar[r] & \dots \ar[r] & D_m^<(\Df_k^>(P)) \ar[r] & 
D_{m-1}^<(\Df_{k+1}^>(P)) \ar[r] & \dots
}
\end{equation}
Differentials of this complex are $A$-homomorphisms and hence its homology, called 
\emph{$\Df$-Spencer homology}, are $A$-modules. They  naturally describe singularities 
of both $A$ (or $\Spec_{\gk}A$) and $P$ (see, for instance, \cite{Boch}).

The \emph{infinite order $\Df$-Spencer complex} $\mathrm{Sp}_*(P)$ 
\begin{equation}\label{d-sp-inf}
\xymatrix@1{
\dots \ar[r] & D_{s}^<(\Df^>P) \ar[r] & D_{s-1}^<(\Df^>(P)) \ar[r] & \dots 
\ar[r] & \Df^>(P) \ar[r] & 0
}
\end{equation}
is the direct limit of natural embeddings 
$\dots\subset \mathrm{Sp_{m}}(P)\subset \mathrm{Sp_{m+1}}(P)\subset\dots$. The homomorphism
$h_*^{\Delta}$ induces a chain map $\gH^{\Delta}\colon \mathrm{Sp}_*(P)\rw\mathrm{Sp}_*(Q)$. 
The related with $\gH^{\Delta}$ homology plays an important role in the geometric theory of PDEs 
(see \cite{Vgeom, KLV}). In particular, in Secondary Calculus computation of this homology for the 
universal linearization operators of (nonlinear) PDEs is the most powerful method for finding higher 
symmetries, conservation laws, field Poisson structures, etc (see \cite{Vper, KLV, KK, SKL}).


\subsection{Higher analogues of multi-derivation functors and $\Df$-Spencer complexes}
Now we can discover many other FDCs just by substituting 
$$
D_{(k)}(P)\df\{\Delta\in \Df_k^<(P)\mid \Delta(1)=0\}
$$ 
for $D(P)$ in the definition of multi-derivation functors. Namely, the
generalization of (\ref{D2}) is as follows:
\begin{equation}\label{hD2}
D_{(k,l)}(P)\df D_{(l)}^<(D_{(k)}(P)\subset \Df_k^>(P))
\end{equation}
\begin{equation}\label{hD2}
\gD_{(k,l)}(P)\df \Df_{l}^<(D_{(k)}(P)\subset \Df_k^>(P))
\end{equation}
and, in general,
\begin{equation}\label{hD2}
D_{(k_1,\dots,k_m)}(P)\df D_{(k_m)}(D_{(k_1,\dots,k_{m-1})}(P)\subset \gD_{k_1,\dots,k_{m-1}}^>(P))
\end{equation}
\begin{equation}\label{DD2}
\gD_{(k_1,\dots,k_m)}^{<>}(P)\df \Df_{(k_m)}^{<>}(D_{(k_1,\dots,k_{m-1})}(P)\subset \gD_{k_1,\dots,k_{m-1}}^>(P))
\end{equation}
A higher analogue of the $\Df$-Spencer complex is associated with a sequence 
$\sigma=(s_1,\dots,s_m), \, s_i\in \Z_+$. It is denoted by $\mathrm{Sp}_{\sigma}(P)$ and 
looks like this:
\begin{eqnarray}
0\rw D_{\sigma_m}(P)\rw\dots\rw D_{\sigma_i}^>(\Df_{m_i}^<(P))\rw 
D_{\sigma_{i-1}}^<(\Df_{m_{i+1}}^>(P))  \rw\dots  \\ 
\dots\rw \Df_{m_0}(P)\stackrel{\rD_{k_0}}{\lrw}P\rw 0  \nonumber
\end{eqnarray}
with $\sigma_i=(s_1,\dots,s_i)$ and $m_i=s_{i+1}+\dots+s_m$. 

If $\sigma\leq\tau$, then  $\mathrm{Sp}_{\sigma}(P)\subset\mathrm{Sp}_{\tau}(P)$. $\Df$-Spencer 
complex can be  also defined for an infinite from the right $\sigma$: 
\begin{equation}\label{}
\dots\rw D_{\sigma_i}^<(\Df^>(P))\rw 
D_{\sigma_{i-1}}^<(\Df^>(P))  \rw\dots\rw \Df(P)\stackrel{\rD_{\infty}}{\lrw}P\rw 0
\end{equation}
For other related FDCs see \cite{V-hom}.


\subsection{Absolute and relative functors}
The previously discussed FDCs may be treated symbolically without reference to concrete algebras
and modules. For instance, by using such symbol as $D^<(D\subset\Df_1^>)$ instead of  
$D^<(D(P)\subset\Df_1^>P)$ we stress that this is an \emph{absolute functor}, i.e., that its construction 
does not depends on the ground algebra $A$. From now on we shall use these symbols. A natural
transformation of one absolute FDC $\Phi$ to another $\Psi$ will be denoted by $\Phi\rw\Psi$.

Natural transformations of the form $\Phi\rw\Psi^<(\Df_k^>)$ are to be distinguished, since they are 
the source of \emph{natural differential operators} (see \cite{MKS}) between representing  $\Phi$
$\Psi$ objects. For example, the naturality property of the exterior differential $d$, namely, that 
$d\circ F^*=F^*\circ d$ for any smooth map $F$, reflects existence of $D_k\rw D_{k-1}^<(\Df_1^>)$
(see section\,\ref{nat}) . 

An example of a \emph{relative functor} of differential calculus is the functor $Q\mapsto \Df_k(P,Q)$ 
for a fixed $A$-module $P$. Relative functors depend on concrete algebras and modules. In
particular, they naturally appear when the ground algebra $A$ or some special modules over 
it have some relevant peculiarities such as Poincare duality, etc.  An instance of that we shall 
see below in connection with \emph{integral forms}. 

It is important to stress that the notion of a FDC include multifunctors. For example,
$\Df_k(\cdot,\cdot)$ is an absolute bifunctor. Multi-functors of the form
$$
\Df_{k_1}^{\varepsilon_1}(P_1, \Df_{k_2}^{\varepsilon_2}(P_2,\dots
\Df_{k_m}^{\varepsilon_m}(P_m,Q)\dots))
$$
with $\varepsilon_i=``<"$ or $``>"$ generalize this simple example. By the space limitations in
this paper we can not give here a due attention to this topics.


\section{The graded generalization}\label{grad}
Naturality of the above general approach also appears in the fact that it automatically generalizes 
to graded commutative algebras and, in particular, to supermanifolds. Indeed, all is needed to this 
end is the graded version of $\delta_a$'s. In particular, this makes algorithmic finding analogues of 
the ``usual" geometrical structures  in the graded context. Moreover, it allows to discover the 
\emph{conceptual meaning} of various well-known quantities, for instance, tensors, that are traditionally defined  by a description. These points are illustrated below. 

\subsection{Differential operators over graded commutative algebras}
Recall that a \emph{graded associative algebra} over a field $\gk$ is a pair $(A, G, \beta)$ where 
\begin{enumerate}
\item $A$ is an associative $\gk$-algebra;
\item $G$ is a commutative monoid written additively; 
\item $A=\oplus_{g\in G}A_g$ with $A_g$'s being $\gk$-vector spaces and 
$A_g\cdot A_h\subset A_{g+h}$;
\end{enumerate}
If $\beta(\cdot, \cdot)$ is an $\dF_2$-valued bi-additive form on $G$, then $A$ is $\beta$-commutative
if $ba=(-1)^{\beta(g,h)}ab$ for $a\in A_g, \,b\in A_h$\footnote{For simplicity we do not consider here a 
more general notion of commutativity (see \cite{MVin, Verb}).}. Similarly, an $A$-module $P$ is \emph{graded} if $P=\oplus_{g\in G}P_g$ and $A_g\cdot P_h\subset P_{g+h}$. In the sequel it will be tacitly assumed that all constructions and operations with graded objects respect gradings. For instance, a
homomorphism $F\colon P\rw Q$ of graded $A$-modules is \emph{graded} of degree $h$ if 
$F(P_g)\subset Q_{g+h}, \,\forall g\in G$. Accordingly, the notation $\Hom_A^h(P,Q)$ will refer to the 
$\gk$-vector space of all $A$-homomorphisms from $P$ to $Q$ of degree $h$ and $\Hom_A(P,Q)$
to the module of all graded $A$-homomorphisms, i.e., $\Hom_A(P,Q)=\oplus_{h\in G}\Hom_A^h(P,Q)$.
An element $p\in P_g$ of a graded module $P$ is called \emph{homogeneous of degree $g$}. We also
shall adopt the simplifying notation $(-1)^{ST}$ for $(-1)^{\beta(\deg S, \deg T)}$ for homogeneous
elements $S$ and $T$ of graded $A$-modules.

If $\Delta\in \Hom_{}\gk(P,Q)$ and $a\in A$ are homogeneous, then we put 
$$
\delta_a(\Delta)=\Delta\circ a_P-(-1)^{a\Delta}a_Q\circ\Delta \quad\mbox{and} \quad
\delta_{a_1,\dots,a_m}=\delta_{a_1}\circ\dots\circ\delta_{a_m}
$$
\begin{defi}\label{GrDO}
Let $P$ and $Q$ be graded modules over a graded commutative algebra $A$. Then
$\Delta\in \Hom_{\gk}(P,Q)$ is a (graded) differential operator of order $\leq k$ if
$$
\delta_{a_0,a_1,\dots,a_k}(\Delta)=0 \quad \mbox{for all homogeneous}  \quad a_0,a_1,\dots,a_k\in A.
$$
\end{defi}

With this definition all above constructions of FDOs, Hamiltonian formalism, etc, automatically 
generalize to the graded case just by literally repeating the corresponding definitions. 


\subsection{Example: Lie algebroids and dioles}
Recall that a \emph{Lie algebroid} over a non-graded commutative $\gk$-algebra $A$ is an 
$A$-module $P$ supplied with a Lie algebra structure $[\cdot,\cdot]_{_P}$ and a homomorphism 
$\alpha\colon P\rw D(A)$, called the \emph{ancor}, such that
\begin{enumerate}
\item $[p,aq]_{_P}=a[p,q]_{_P}+\alpha(p)q, \,\forall a\in A, \,p,q\in P;$
\item $\alpha$ is a Lie algebra homomorphism from $(P, [\cdot,\cdot]_P)$ to $(D(A), [\cdot,\cdot])$.
\end{enumerate}
This is a general algebraic version of the standard definition when $A=C^{\infty}(M)$ and 
$P=\Gamma(\pi)$ with $\pi$ being a vector bundle over $M$.

This at the first glance not very usual  geometrical object is, in fact, a graded analogue of a Poisson 
manifold. More exactly, the corresponding $\Z$-graded algebra $\cA$, called the algebra of 
\emph{diols} or \emph{diole algebra},  is defined to be
\begin{enumerate}
\item $\cA_0=A, \,\cA_1=P$ and $\cA_i=\{0\}, \,i\neq 0,\,1$;
\item the product in $\cA_0\subset\cA$ is that in $A$ and the product 
$\cA_0\cdot\cA_1\subset \cA_1$ is the $A$-module product $(A,P)\rw P$.
\end{enumerate}
Note that $\cA_1\cdot\cA_1\subset\cA_2=\{0\}$ and that $\cA$ is graded commutative with respect 
to the trivial sign form $\beta$. See \cite{Misa} for further details.

Let  $\{\cdot,\cdot\}$ be a Poisson structure in $\cA$ of degree $-1$, i.e., a graded Lie algebra 
structure in $\cA$ such that $\{P,P\}\subset P, \,\{P,A\}\subset A$, which is additionally a biderivation 
of $\cA$.  The ``Hamiltonian" field $\alpha(p)\df\{p,\cdot\}|_A, \,p\in P,$ is, obviously, a derivation of 
$A$. From the biderivation property $\{ap,b\}=\{a,b\}p+a\{p,b\}, \, a,b\in A, \,p\in P,$ and 
$\{A,A\}\subset\cA_{-1}=\{0\}$ we see that $\alpha(ap)=a\alpha(p)$, i.e., that $\alpha\colon P\rw A$ 
is a homomorphism of $A$-modules. Then the Jacobi identity
$$
\{\{p,q\},a\}+\{\{a,p\},q\}+\{\{q,a\},p\}=0 \quad\Leftrightarrow \quad\alpha(\{p,q\})(a)=[\alpha(a),\alpha(q)](a)
$$
tells that $\alpha$ is a homomorphism of Lie algebras. Similarly, by putting 
$[\cdot,\cdot]_{_P}\df \{\cdot,\cdot\}|{_P}$ we see that condition (1) in the definition of algebroid is 
exactly the Leibniz rule $\{p,aq\}=\{p,a\}q+a\{p,q\}$. So, Lie algebroids are nothing but Poisson 
structures on algebras of dioles of degree $-1$ or \emph{fat Poisson manifolds} in the sense of
\cite{fat}. An advantage of this interpretation is that it puts algebroids into the rich context of 
differential calculus over diole algebras and, in particular, makes obvious the analogy with the 
standard Poisson geometry. 

In this connection one may be curious about Poisson structures of different degrees over $\cA$. 
Among these only structures of degrees from $-2$ to $1$ nay be nontrivial. Their ``non-graded" 
description is as follows.

A Poisson structure of degree $-2$ is just an $A$-bilinear and $A$-valued skew-symmetric  form on 
$P$. Poisson structures of degree $1$ are elements of $D_2(P)$. Structures of degree $0$ are more
 complicated.  Each of them consists of a Poisson structure $\{\cdot,\cdot\}_{_A}$ on $A$ and a flat 
 \emph{Hamiltonian connection}, which lifts the ``Hamiltonian vector field" 
 $X_a\df \{a, \cdot\}_{_A}, \,a\in A,$ to the derivation $\nabla_a\colon P\rw P$ of $P$ over $X_a$. This 
 means that $\nabla_a(bp)=X_a(b)p+b\nabla_a(p), \,b\in A, \,p\in P$. Additionally, it is required that
 $\nabla_{ab}=a\nabla_b+b\nabla_a$.


\section{Representing objects}\label{rep}
In our approach covariant tensors, jets and other covariant objects of the standard differential geometry appear as elements of objects representing FDCs in suitable subcategories of the category 
$A\mathbf{Mod}$ of $A$-modules. Following are main details of this construction.


\subsection{Representing objects$\colon$generalities}\label{gen-ob}
For simplicity we shall consider only non-graded case.
In the first approximation an object representing a FDC $\Phi$ in a category $\cK$ of  $A$ modules
is an $A$-module $\cO_{\cK}(\Phi)$ such that $\Phi$ is equivalent to the functor
$P\mapsto \Hom_A(\cO_{\cK}(\Phi),P), \,P\in\Ob\cK$. Representing the same functor objects are
naturally isomorphic. The homomorphism $\cO_{\cK}(\Phi)\rw P$ representing $\theta\in\Phi(P)$
will be denoted by $h_{\theta}$. A natural transformation $\Phi\rw \Psi$ of FDCs is then represented 
by a homomorphism $\cO_{\cK}(\Psi)\stackrel{\Upsilon}{\lrw}\cO_{\cK}(\Phi)$ of $A$-modules. Namely, 
by identifying
$\Phi(P)$ and $\Hom_A(\cO_{\cK}(\Phi),P)$, etc, we have:
$$
\Hom_A(\cO_{\cK}(\Phi),P)\lrw\Hom_A(\cO_{\cK}(\Psi),P), \quad h_{\theta}\mapsto 
h_{\theta}\circ\Upsilon
$$

Let $A$ and $B$ be commutative algebra. For an $A$-module $P$, an $B$-module R and an 
$(A,B)$-bimodule $Q$ we have the canonical isomorphism 
\begin{equation}\label{tens-hom}
\Hom_B(P\otimes_AQ, R)=\Hom_A(P, \Hom_B(Q, R)).
\end{equation}
An obvious consequence of it is that $\cO_{\cK}(\Phi)\otimes_A\cO_{\cK}(\Phi)$ represents the
composition $\Phi\circ\Psi$, i.e., the functor $P\mapsto \Phi(\Psi(P))$. 

Assume now that $\Psi$ is $A$-bimodule-valued. Then by labeling the corresponding two 
$A$-module structures by $``<"$ and $``>"$ we obtain two $A$-module-valued functors, 
$\Psi^<$ and $\Psi^>$. Each of these multiplication by $a\in A$ is a natural transformation 
of $\Psi$ and hence induce an endomorphism of $\cO_{\cK}(\Psi)$. This way $\cO_{\cK}(\Psi)$
acquires an $A$-bimodule structure. Accordingly, the corresponding  $A$-module structures 
will be also denoted by $``<"$ and $``>"$. The $A$-module $\Phi^<(\Psi^>(P))$ is defined as 
the $\gk$-vector space coinciding with $\Phi(\Psi^>(P))$, in which the $A$-module structure is 
induced by the $<$-module structure in $\Psi(P)$. This way we get the functor $\Phi^<(\Psi^>)$.
Once again it follows from isomorphism (\ref{tens-hom}) that the representing object for
$\Phi^<(\Psi^>)$ is $^<\cO_{\cK}(\Psi)_{>}\otimes_{A}\cO_{\cK}(\Phi)$. This symbol tells that
the tensoring is taken with respect to the $>$-structure of $\cO_{\cK}(\Psi)$, while the 
$A$-module structure of the obtained tensor product is induced by the $<$-structure in 
$\cO_{\cK}(\Psi)$.


\subsection{Existence of representing objects}\label{ex-rep}
Representing objects in the category $A\mathbf{Mod}$ of all $A$-modules exist for all FDCs. The 
techniques sketched in the previous subsection reduce, basically, their construction to that for
functors $\Df_k(P,\cdot)$ (see \cite{KLV, Vez}). These objects are called \emph{$k$-th order jets} 
of $P$ and are denoted by $\cJ^k(P)$. The construction of these $A$-modules is rather elementary.
Indeed, consider with this purpose the $A$-module $A\otimes_{\gk}P$ and associate with an 
$a\in A$ the homomorphism
$$
\delta^a\colon A\otimes_{\gk}P\rw A\otimes_{\gk}P, \quad 
\delta^a(a'\otimes_{\gk}p)=a'\otimes_{\gk}ap-aa'\otimes_{\gk}p.
$$
Next, denote by $\mu_{k+1}$ the submodule of  $A\otimes_{\gk}P$ generated by all elements of
the form $(\delta^{a_0}\circ\delta^{a_1}\circ\dots\circ\delta^{a_k})(a\otimes_{\gk}p)$ and put 
$$
\cJ^k(P)\df \frac{A\otimes_{\gk}}{\mu_{k+1}} \quad\mbox{and} \quad j_k\colon P\rw \cJ^k(P), 
\quad j_k(p)=1\otimes_{\gk}p \mod \mu_{k+1}.
$$
Then we have
\begin{proc}\label{jet}
For any $\Delta\in\Df_k(P,Q)$ there is an unique $A$-homomorphism $h^{\Delta}$ that makes
the following diagram commutative:
\begin{equation}\label{j-k}
\xymatrix{
P \ar[r]^-{j_k} \ar[rd]_{\Delta} & \cJ^k(P) \ar[d]^{h_{\Delta}} \\
& Q
}
\end{equation}
The correspondence $\Delta \mapsto h_{\Delta}$ establishes an isomorphism of $A$-modules 
$\Df_k^<(P,Q)$ and $\Hom_A(\cJ^k(P),Q)$. Moreover, $\cJ^k(P)=^<\cJ^k(A)_>\otimes_AP$.
\end{proc}
As it is easy to see, $j_k$ is a $k$-th order DO and  proposition\,\ref{jet} tells that it is \emph{universal}.
According to subsection\,\ref{gen-ob}, the embedding of functors 
$\Df_l(P,\cdot)\rw\Df_k(P,\cdot), \,l\leq k,$ induces a natural projection 
$\pi_{k,l}\colon\cJ^k(P)\rw\cJ^l(P)$.

 Conceptually,  \emph{$k$-th order differental forms} are defined to be elements of the $A$-module representing functor $D_k$ in $A\mathbf{Mod}$, which is denoted by $\Lambda^k(A)$. This module
 may be constructed by the methodas  above. A natural splitting
 $\Df_1^<=\id\oplus D$ where $\id=\Df_0^<$ is the identity functor suggests to define 
 $\Lambda^1(A)$ to be the kernel of $\pi_{1,0}\colon\cJ^1(A)\rw\cJ^0(A)=A$. In this approach,
 the differential $d\colon A\rw \Lambda^1(A)$ is defined by the formula $da\df j_1(a)-aj_1(1_A)$.
 
 It follows from the last construction in subsection\,\ref{gen-ob} and proposition\,\ref{jet} that
 the representing functor $D_{k-1}^<(\Df_1^>)$ $A$-module is $\cJ^1(\Lambda^{k-1}(A))$. 
 The analogous to (\ref{j-k}) diagram
 $$
\xymatrix{
\Lambda^{k-1}(A) \ar[r]^-{j_k} \ar[rd]_{d} & \cJ^k(\Lambda^{k-1}(A)) \ar[d]^{h^{}} \\
& \Lambda^{k}(A),
}
$$
in which $h$ is the homomorphism representing the transformation of functors 
$D_k\rw D_{k-1}^<(\Df_1^>)$, is the definition of the exterior differential $d$. Now it is not difficult 
to see that the \emph{jet-Spencer complex}
\begin{equation}\label{j-sp}
0\lw \Lambda^n(A)\lw\dots \cJ^k(\Lambda^{n-k}(A))\lw\cJ^{k+1}(\Lambda^{n-k-1}(A))\lw\dots
\end{equation}
represents the functor $P\mapsto\mathrm{Sp_{n}}(P)$ (see \ref{d-sp}).  For more examples of this 
kind see \cite{KLV} (chapter I) and \cite{V-hom, Vez, V}.

It is possible to construct an $A$-module  representing a  single FDC $\Phi$ in a category $\cK$
as the quotient module $\cO_{A\mathbf{Mod}}(\Phi)/K$ where $K$ is the intersection of kernels
of all homomorphisms $\cO_{A\mathbf{Mod}}(\Phi)\rw Q, \,Q\in\Ob\cK$. However, this module
does not, generally, belong to $\cK$. This makes impossible to represent in $\cK$ all FDCs and
connecting them natural DOs. Categories of $A$-modules that contain so-defined single
representing $A$-modules are called \emph{differentially closed} (see \cite{Vez} for more details).

An important example of differentially closed categories is the category $\Gamma A\mathbf{Mod}$ 
of geometric $A$-modules over a geometrical algebra $A$. Representing $A$-modules in this category 
are  geometrizations of representing $A$-modules in $A\mathbf{Mod}$ (see subsection\,\ref{gm}).
Importance of  the category of geometric modules is that it is in full compliance with
the observability principle. Indeed, we have
\begin{proc}\label{rpr}
Representing objects in the category of geometric $A$-modules over $A=C^{\infty}(M)$ are 
identical to the corresponding objects in the standard differential geometry.
\end{proc}
This means that differential forms, jets, etc, in the ordinary sense of these terms are nothing else
than elements of $C^{\infty}(M)$-modules that represents the corresponding FDC in the category of
geometric modules. The following example illustrate the drastic difference between categories of
all and geometric $C^{\infty}(M)$-modules.
\begin{ex}
Let $d_{alg}$ denote the exterior differential in $A\mathbf{Mod}$.  If $A=C^{\infty}(\R)$, then
$d_{alg}(e^x)\neq e^xd_{alg}x$. In other words, 
$(d_{alg}(e^x)-e^xd_{alg}x)\in \mathrm{Ghost}(C^{\infty}(\R))$.
\end{ex}
The reader will find more about in \cite{MVhost}.


\subsection{Naturality of $d$ and related topics}\label{nat}
A remarkable property of ordinary differential forms, jets, etc, is their \emph{naturality}. This means
that any smooth map $F\colon M\rw N$ is accompanied by a map 
$F^*\colon \Lambda^i(N)\rw \Lambda^i(M)$. We denote by $\Lambda^i(L)$ the $C^{\infty}(L)$-module
of (ordinary) $i$-th order differential forms on the manifold $L$. This ``experimental" fact has the 
following explanation.

Let $H\colon A_1\rw A_2$ be a homomorphism of commutative algebras and $\Phi$ an absolut FDC.
In this situation any $A_2$-module $Q$ acquires an $A_1$-module structure with the $A_1$-module multiplication $(a_1,q)\mapsto H(a_1)q, \,a_1\in A_1, \,q\in Q$. 

Now assume  that $\cK_i, \,i=1,2,$ is a differentially closed category of $A_i$-modules and that any 
$Q\in \Ob\cK_2$ belongs to $\cK_1$ as $A_1$-module. Also change the notation by putting  
$\Lambda_{\cK}^i(A)=\cO_{\cK}(D_i)$, $\cJ^k_{\cK}(A)=\cO_{\cK}(\Df_k^<)$ and denote by 
$d_{\cK}\colon \Lambda_{\cK}^i(A) \rw \Lambda_{\cK}^{i+1}(A)$ the exterior differential in $\cK$. 
The composition $X=d_{\cK_2}\circ H$ is a derivation of $A_1$ with values in 
$\Lambda_{\cK_2}^1(A_2)$ considered as an $A_1$-module in $\cK_1$. By universality 
of $d_{\cK_1}, \;X=h_X\circ d_{\cK_1}$. So, by putting $H_{\Lambda^1}\df h_X$ we get 
a commutative diagram at the left :
$$
\xymatrix{
\Lambda_{\cK_1}^1(A_1) \ar[r]^-{H_{\Lambda^1}} &\Lambda_{\cK_2}^1(A_2) &\ar @{} [dr] |{\Longrightarrow}&& \Lambda^1(N) \ar[r]^-{F^*} & \Lambda^1(M)\\
 A_1 \ar[r]^-{H} \ar[u]^{d_{\cK_1}} & A_2 \ar[u]_{d_{\cK_2}} & &
 &C^{\infty}(N) \ar[r]^-{F^*} \ar[u]^-{d}& C^{\infty}(M) \ar[u]_-{d}
}
$$
The diagram at the right expressing naturality of 1st order differential forms and the exterior
differential $d$ is the specialization of the left diagram to $H=F^*$ and  categories of geometrical 
modules  $\cK_i$'s in view of theorem\,\ref{spc} and proposition\,\ref{rpr}.

The same arguments explain naturality of modules $\cJ^k_{\cK}(A)$ and the DO 
$j_k\colon A\rw \cJ^k_{\cK}(A)$. Also, together with inductive arguments used in the definition of
functors $D_k$s they explain naturality of higher order differential forms and exterior differentials,
jet-Spencer complexes and so on. Higher analogues of de Rham and Spencer complexes are
examples of natural differential operators, which can be hardly  discovered by traditional methods
(compare with \cite{MKS}). 


\subsection{Multiplicative structure in $\cJ^k$}
Representing modules may have various additional structures coming from specific natural 
transformations of FDCs. This point is illustrated below. Assuming that a differentially closed category
$\cK$ is fixed we shall omit the subscript $\cK$ in the notation of representing modules in this category.

From the defining $D_{k+1}$  formula $D_{k+1}=D^<(D_{k}\subset \Df_1^>)$ follows the inclusion
transformation $D_{k+1}\rw D\circ D_k$. By iterating this procedure we can construct an inclusion 
$D_{k+l}\rw D_l\circ D_k$. Then, according to the general principles of subsection\,\ref{gen-ob},
the last inclusion is represented by a homomorphism 
$\wedge_{k,l}\colon \Lambda^k(A)\otimes_A \Lambda^l(A)\rw \Lambda^{l+l}(A)$ of representing
modules. If $\cK$ is the category of geometric modules over $C^{\infty}(M)$, then $\wedge_{k,l}$
is the standard wedge product of differential forms on $M$. As before, the naturality of the 
so-defined wedge product straightforwardly follows from its definition.

Our another example concerns modules of jets. The \emph{diagonal transformation} of 
functors $\Df_k^<\stackrel{\gm_k}{\lrw} \Df_k^<\circ \Df_k^<$:
\begin{equation}\label{}
 \Df_k^<P\ni \Delta \mapsto\gm_k(\Delta)\in \Df_k^<(\Df_k^<P), \quad\gm_k(\Delta)(a)=\Delta\circ a_P
\end{equation}
is represented by the homomorphism $\cJ^k(A)\otimes_A\cJ^k(P)\stackrel{\gm^k}{\lrw}\cJ^k(P)$
of the representing modules. This supplies $\cJ^k(P)$ with a natural $\cJ^k(A)$-module structure.
It is easily deduced from this definition that
$$
j_k(a)\cdot j_k(p)=j_k(ap), \;a\in A, p\in P  \quad\mbox{where}  
\quad j_k(a)\cdot j_k(p)\df\gm^k(j_k(a)\otimes j_k(p))
$$
In particular, if $P=A$, then $\gm^k$ supplies $\cJ^k(A)$ with a commutative algebra  structure.

The standard operator of insertion of a vector field $X$ to differential forms is due to a natural 
transformation of \emph{relative} functors $i^X\colon D_k\rw D_{k+1}, \,X\in D(A)$. For instance, 
the insertion operator $i_X\colon \Lambda^2\rw\Lambda^1$ represents the transformation
of functors
$$
i^X\colon D(P)\ni Y\mapsto \{a\mapsto X(a)Y-Y(a)X\}\in D_2(P).
$$
Similarly can be defined the Lie derivative of differential forms as well as higher order analogues of it
and of the insertion operation (see \cite{Vezo, Vez}).


\subsection{Tensors conceptually} According to the standard coordinate-free definition, covariant 
tensors, are $C^{\infty}(M)$-multilinear functions on vector fields on $M$.  This definition, being 
descriptive, does not tells anything about the role of these objects in the structure of differential 
calculus. For instance, it does not explain why skew-symmetric tensors, i.e., differential forms, 
are related by a natural DO, the exterior differential, while the symmetric tensors do not. 
Another similar question is why a natural connection, namely, that of  Levy-Civita, is associated 
only with non-degenerate symmetric tensors. In this regard it is instructive noticing that  various  
attempts to construct an analogue of the Levy-Civita connection for symplectic manifolds have 
been  made not long ago without any positive result.

Below 
we shall sketch the conceptual approach to tensors and, as a byproduct, shall answer the above 
two questions. The decisive idea is to interpret covariant tensors as special first order differential 
forms over the algebra of \emph{iterated differential forms} (see \cite{iter1}). For simplicity we shall 
discuss only the non-graded case.

We begin from the algebra $\Lambda_1\df\Lambda_{\cK}^*(A)$ of differential forms in a suitable 
differentially closed category of $A$-modules, for instance, the category of geometrical modules 
over $C^{\infty}(M)$. In order to save ``space-time" we shall omit all references to $\cK$ including the 
notation. $\Lambda_1$ is a $\Z$-graded commutative algebra  with the ordinary multiplication 
$\Z\times\Z\rw\Z$ for the grading form. Therefore, the algebra of differential forms over $\Lambda_1$
denoted $\Lambda_2\df\Lambda^*(\Lambda_1)$ is well-defined. The exterior differential $d=d_1$
is a derivation or a ``vector field" over $\Lambda_1$. It naturally extends as the Lie derivative $L_{d_1}$
to $\Lambda_2$. Hence $L_{d_1}$ commutes with the exterior differential $d_2$ in $\Lambda_2$.
It is convenient to abuse the notation by using $d_1$ for $L_{d_1}$. Note also that 
$\Lambda_1=\Lambda^0(\Lambda_1)\subset\Lambda_2$. By continuing this process we inductively
construct the algebra $\Lambda_k\df\Lambda^*(\Lambda_{k-1})$ of \emph{$k$-times iterated 
differential forms} with commuting differentials $d_1,\dots, d_k$. Natural inclusions 
$\Lambda_{k-1}=\Lambda^0(\Lambda_{k-1})\subset\Lambda_k$ allow to define the \emph{algebra
of iterated differential forms} 
$$
\Lambda_{\infty}\df\bigcup_{0\leq k}\Lambda_k \quad\mbox{with} \quad \Lambda_0=A.
$$ 
This algebra is \emph{conceptually closed}, i.e., $\Lambda^*(\Lambda_{\infty})$ is naturally isomorphic
to $\Lambda_{\infty}$.

In the category of geometric modules over $A=C^{\infty}(M)$ the above-mentioned interpretation of
covariant tensors as iterated differential forms looks very simply and is as follows:
$$
i_p \colon T^k(M)\ni df_1 \otimes_A df_2 \otimes_A\dots\otimes_A df_k \quad\mapsto \quad
d_1 f_1\wedge d_2 f_2\wedge\dots\wedge d_k f_k \in \Lambda_k
$$
In particular, if $\tau'=\tau_{ij}dx^i\otimes dx^j\in T^2(M)$, then 
$\tau\df i_2(\tau')=\tau_{ij}d_1x^id_2x^j\in \Lambda_2^1$ and
$$
d_2(d_1\tau)=\tau_{ij,kl}d_1x^id_1x^kd_2x^jd_2x^l - \gamma_{ijk}d_1x^id_2x^jd_2(d_1x^k)-
g_{ij}d_2(d_1x^i)d_2(d_1x^j)
$$
where $\tau_{ij,m}=\p\tau_{ij}/\p x^m$, etc,  $\gamma_{ijk}=\tau_{ik,j}+\tau_{kj,i}-\gamma_{ij,k}$,
$g_{ij}=1/2(\tau_{ij}+\tau_{ji})$ and we omit symbols of wedge products. Note that $\gamma_{ijk}$ 
is the double fist kind Christoffel symbol of the Levi-Civita connection if the tensor $\tau'$ is 
symmetric, i.e., $\tau'=g\df g_{ij}dx^i\otimes dx^j$.

Assume now that $g$ is non-degenerate. Then the map $\mathrm{grad_g}\colon f\mapsto \mathrm{grad_g}f$
is well-defined. It is a $D(A)$-valued derivation of $A$ and hence $\mathrm{grad_g}=h_g\circ d$
with $h_g=h^{\mathrm{grad_g}}\colon \Lambda^1(A)\rw D(A)$. The isomorphism $h_g$ naturally 
extends to an $A$-linear derivation of $\Lambda_1=\Lambda^*(A)$ with values in the $\Lambda^*(A)$-module $\Lambda^*(A)\otimes_AD(A)$ denoted $h_g^*$. Similarly, $h_g^*$, being a derivation of 
$\Lambda_1$ induces a $\Lambda_1$-linear derivation 
$$
h_g^{(2)}\colon\Lambda_2\rw \Lambda_2\otimes_{\Lambda_1}(\Lambda_1\otimes_AD(A))=
\Lambda_2\otimes_AD(A).
$$
In coordinates $h_g^{(2)}=g^{l\alpha}i_{\p_{\alpha}}\otimes_A i^{(2)}_{\p_{\alpha}}$ with
$i^{(2)}_{\p_{\alpha}}=i_{i_{\p_{\alpha}}}$. Then
\begin{equation}\label{L-C}
\Gamma(\tau)\df=\frac{1}{2}h^{(2)}_{g}(d_2(d_1\tau))=
(d_2(d_1x^{\alpha})+\Gamma^{\alpha}_{ij}d_1x^id_2x^j)\otimes_Ai_{\p_{\alpha}}
\end{equation}
where $\Gamma^{\alpha}_{ij}=1/2g^{\alpha k}\gamma_{kji}$. If $\tau'=g$, then
$\Gamma^{\alpha}_{ij}$'s are the Christoffel symbols of the pseudo-metric $g$.
Call $\Gamma(\tau)$ the \emph{Levi-Civita form of $\tau$}. Since $\Gamma(\tau)$ is a
vector-valued graded form, its graded Frolicher-Nijenhuis square $[\cdot,\cdot]^{FN}$ is well-defined
and gives the \emph{curvature form}:
\begin{equation}\label{curv}
[\Gamma(\tau),\Gamma(\tau)]^{FN}=R^{\alpha}_{ijk}d_1x^kd_2x^jd_2x^i\otimes_Ai_{\r_{\alpha}}
\end{equation}
with $R^{\alpha}_{ijk}=\p_{i}\Gamma^{\alpha}_{jk}-\p_{j}\Gamma^{\alpha}_{ik}+
\Gamma^{\alpha}_{i\beta}\Gamma^{\beta}_{jk}-\Gamma^{\alpha}_{j\beta}\Gamma^{\beta}_{ik}$.
Similarly, all other standard quantities related with the Levi-Civita connection can be obtained by
applying to $\Gamma(\tau)$ natural operators of differential calculus over $\Lambda_1$ (see 
\cite{iter2}). So, all these facts lead to recognize that, conceptually, $\Gamma(\tau)$ is what 
should be called the Levi-Civita connection associated with a second order covariant tensor 
with non-degenerate symmetric part. This interpretation inserts Riemannian geometry into the 
rich machinery of differential calculus over iterated differential forms and, in particular, allows to
avoid seemingly natural questions related with tensors, which are, in fact, conceptually ill-posed 
as, for example, the question about analogue of the Levi-Civita connection for symplectic manifolds.


\subsection{An example of application: natural equations in general relativity}
A richer and conceptually certain mathematical language offers more possibilities to mathematically
formalize various situations in physics. This is especially important when the subject is of a non-intuitive
character. We shall illustrate this common place with an application of iterated forms to general relativity
(see \cite{iter2}). 

Let $\tau=g+\omega$ be a covariant second order tensor field
on a 4-fold $M$ with $g$ and $\omega$ being its symmetric and skew-symmetric parts, respectively.
We shall interpret   $\tau, g$ and $\omega$ as iterated $(1,1)$-forms. In the context of general relativity
it is natural to interpret $g$ as the metric in the space shaped by the matter field $\omega$,
i.e., by the ``fermionic" part of $\tau$. By analogy with the classical vacuum Einstein equation
$\mathrm{Ric}(g)=0$ we assume that $g$ and $\omega$ are connected by the equation 
$\mathrm{Ric}(\tau)=0$ where $\mathrm{Ric}(\tau)$ is the Ricci tensor of the connection 
$\Gamma(\tau)$. Describe it in coordinates.

If $R^{\alpha}_{ijk}$'s are as in (\ref{curv}), then $\mathrm{Ric}(\tau)=R_{ik}d_1x^id_2x^k$ with
$R_{ik}=R^j_{ijk}$. Denote by $\nabla_g$ the covariant differential of the connection $\Gamma(g)$
and put  $\mathrm{Ric}(g)=R_{ik}^{(g)}d_1x^id_2x^k$. Then $\mathrm{Ric}(\tau)=0$ reads as
follows:
\begin{eqnarray}\label{cosmo}
R_{ij}^{(g)}+\frac{9}{16}g^{kl}g^{mn}\p_{[m}\omega_{il]}\p_{[k}\omega_{jn]}=0 \\
\nabla_g(\p_{[i}\omega_{jk]})=0 \qquad\qquad\nonumber
\end{eqnarray}
The second of these equations describes a perfect irrotational fluid. A remarkable feature of known
exact solutions of this equation is that they describe an expanding universe. At the moment physical interpretation of this fluid is not very clear. Speculatively,  one may think that ``molecules" forming it 
are galaxies.  By concluding we stress that $\omega$ is a rather simplified model of matter. But, on 
the other hand, various richer models can be proposed by varying the algebra of observables. 


\subsection{Integration}\label{int}
In the category of smooth orientable manifolds \emph{integrands}, i.e., the quantities to be integrated, 
are top differential forms. In the case of supermanifolds differential forms can be of any positive degree 
and, therefore, none of them can not be considered as integrand. By introducing \emph{integral forms}
F. Berezin had overcome this difficulty (see \cite{Ber}). Further development of the original Berezin
approach have led to a general construction of integral forms over general graded commutative 
algebras (see \cite{Verb}). The idea of this construction is to pass  from the complex of differential
forms
\begin{equation}\label{dRm}
0\rw A=\Lambda^0(A)\stackrel{d_0}{\lrw} \Lambda^1(A)\stackrel{d_1}{\lrw}\dots
\stackrel{d_1}{\lrw}\Lambda^i(A)\stackrel{d_1}{\lrw}\dots
\end{equation}
to the complex of adjoint operators 
\begin{equation}\label{int-form}
0\lw \Sigma^0(A)\stackrel{\widehat{d_0}}{\llw} \Sigma^1(A)\stackrel{\widehat{d_1}}{\lw}\dots
\stackrel{\widehat{d_i}{}_{-1}}{\llw}\Sigma^i(A)\stackrel{\widehat{d_i}}{\llw}\dots
\end{equation}
If $A=C^{\infty}(M), \,\dim\,M=n$ and $\hat{P}\df\Hom_A(P, \Lambda^n(A))$, then 
$\Sigma^i(A)=\hat{\Lambda}^i(A)=\Lambda^{n-i}(A)$ and $\widehat{d_i}=\pm d_{n-i-1} $. In particular,
$\Sigma^0(A)=\Lambda^n(A)$ (see \cite{VL, VCSP}) and elements of $\Sigma^0(A)$ may be 
interpreted as integrands. This is crucial, since these elements belong to the initial part of a complex
and not to its top term, which does not, generally, exist. Not very satisfactory point here is that
the definition of  $\hat{P}$ is based on existence the top term $\Lambda^n(A)$. Fortunately,
this inconvenience can be resolved by observing that $\hat{P}$ is the cohomology of the complex
$\Df(P,\Lambda(A))$ of $A$-homomorphosms
\begin{eqnarray}\label{dc}
0\rw \Df^>(P,A)\stackrel{w_P}{\lrw}\Df^>(P,\Lambda^1(A))\stackrel{w_P}{\lrw}\dots
\stackrel{w_P}{\lrw}\Df^>(P,\Lambda^i(A))\stackrel{w_P}{\lrw}\dots, \\
w_P(\Delta)\df d\circ\Delta,   \quad  \Delta\in\Df(P,\Lambda(A)) \qquad\qquad\quad\nonumber
\end{eqnarray}
denoted by $\hat{P}=H(w_P), \,\hat{P}=\oplus_i\hat{P}^i$  with $\hat{P}^i=H^i(w_P)$. Note that
complex (\ref{dc}) is well-defined for any $G$-graded algebra and that $\hat{P}$ is 
$(\Z\times G)$-graded. If $A=C^{\infty}(M)$, then the only nontrivial component of $\hat{P}$ is
$\hat{P}^n$. The module $\hat{P}$ is called the \emph{adjoint} to $P$.

A DO $\square\in\Df(P,Q)$ induces a cochain map $\Df(Q,\Lambda(A))\rw\Df(P,\Lambda(A))$
and the induced map in cohomology $\hat{\square}\colon\hat{Q}\rw\hat{P}$ is called the 
\emph{adjoint} to $\square$ operator.   With these definitions the adjoint to the de Rham complex
(\ref{int-form}) is well-defined with $\Sigma^i(A)\df\widehat{\Lambda^i(A)}$. Elements of $\Sigma^s(A)$
are called \emph{$s$-order integral forms}. The $A$-module $\mathcal{B}(A)\df\Sigma^0(A)$ is 
called the \emph{Berezinian} (in the category $\cK$ of $A$-modules) (see \cite{Verb}). Originally,
integral forms and Berezinians appeared in the context of supermanifolds (see \cite{Ber, Vor, Pen}) 
and this is the generalisation of these notions to arbitrary graded  commutative algebras. 
\begin{rmk}
For some ``good" algebras $A$ the Berezinian can be thought as an 1-dimensional projective 
$A$-module supplied with a flat right connection $($see \cite{Man, MMV}$)$. This simple interpretation
could be convenient in practical computations but is not satisfactory as a conception.
\end{rmk}
Computations of Berezinians for many algebras of interest essentially proceed along the same lines 
as in \cite{VL, VCSP}. For example, with this method can be described the Berezinian of the algebra 
iterated forms $\Lambda_k(A)$ for $A=C^{\infty}(M)$. It turns out that the only nontrivial homogeneous
component of $\mathcal{B}(\Lambda_k)$ is of order $2^{k-1}n$, which is canonically isomorphic to  
$\Lambda_k$ (see \cite{iter3}).

Finally, in this approach integration is the map that associates with an integral form 
$\omega\in\mathcal{B}(A)$ its homology class in complex $(\ref{int-form})$.


\section{Conclusions}
In this quick trip through differential calculus over commutative algebras we have tried to show 
expressive capacity and universality of this new language. The based on it standpoint forces  
some different views on both traditional and in development parts of mathematics and theoretical 
physics. We shall outline some of them by starting from the following historical parallel with the 
purpose to avoid a formal theorizing.
\subsection{Sheaves as ``broken lines"}
Since there are no means to study arbitrary (smooth) curves and other ``curvilinear objects" within the
Euclidean geometry framework, geometers in antiquity reached a psychological comfort with the idea
that a curve is the limit of  inscribed in it broken lines. Even if the intuitive term ``limit" have been well 
defined this would not be a self-consistent definition. Nevertheless, this intuition helped to compute
the length or other geometrical characteristics of some particular curves before the invention of 
differential calculus. Moreover, broken-lines-like considerations were among decisive factors that
had led to discovery of differential calculus.

This is one of many situations one may meet either in the history or in contemporary mathematics 
when mathematical objects are constructed from still treatable in the old language pieces while the 
new adequate language is not yet in hand. Sheaves are well-known examples of this kind. For 
instance, the necessity of sheaves in modern complex geometry is due to the fact that a complex 
manifold is defined as something sewn from open pieces of $\C^n$. As a consequence all other
relevant geometrical quantities in complex geometry are defined as cohomology of suitable sheaves.
The passage to limit in the definition of sheaf cohomology is pretty parallel to the curves understood
as limits of broken lines. Implicitly, this reflects the fact that a complex manifold is not, generally, the
spectrum of its algebra of holomorphic functions.  In this sense complex manifolds are not ``observable".
But the observability can be immediately reached if a complex manifold is defined to be a smooth manifold 
supplied with an integrable Nijenhuis tensor. In this approach all necessary holomorphic objects, 
for instance, tensors, are defined as compatible with the Nijenhuis tensor ones. This way sheaves 
can be eliminated  from complex geometry together with the corresponding heavy technical 
instruments like derived functors, etc. This not only simplifies and enriches the theory but also 
leads to new important generalizations. As an example we mention the theory of singularities of 
solution of PDEs where analogues of complex geometry are of fundamental importance  
(see \cite{Vmonge}).

This parallel between sheaves and broken lines allows to foresee that the language 
of differential calculus over commutative algebras will inevitably substitute the language 
of sheaves  in the future . At the same time the fundamental role of ``sheaf technologies" 
in the past should be highly recognised. In particular, they implicitly contributed to preparing 
the land for differential calculus over commutative algebras.


\subsection{Some general expectations}
Essentially the same arguments as in the preceding subsection are applied as well to other areas 
of contemporary mathematics and theoretical physics. We shall briefly indicate those of them where
all-round implementation of the DCCA-based methods looks most promising. \newline

- {\it Algebraic geometry}.  Algebraic geometry is an area, which suggests itself introduction of  
DCCA. For that is sufficient to change the object of study by passing, according to the ``philosophy 
of observability", from affine or projective varieties to the corresponding algebras. This allows  direct application of methods of differential geometry. For instance, the Spencer cohomology of algebraic
varieties, defined as in subsections \ref{Spen} and \ref{ex-rep}, are fine invariants of their singularities
(see \cite{Boch})/ As such it suggests an alternative approach to the problem of resolution of singularities.
Also, existence of singularities does not prevent direct definition of the De Rham and other basic 
cohomology of algebraic varieties. For some other examples see \cite{acta}. This is shortly why it 
is natural to think that a systematic review of foundations of algebraic geometry on the basis of DCCA 
would give a strong new impetus to this classical area. \newline

- {\it Geometry of PDEs and secondary calculus}. Originally, one of the main stimuli to develop DCCA
sprang out of the necessity to construct the complete analogue of differential geometry on the ``space
of all solutions of a given PDE" (see \cite{KLV}). The intuitive idea of such a``space" is formalised with 
the concept of a \emph{diffiety}.  Diffieties form a special class of, generally, infinite-dimensional 
manifolds where traditional methods and means fail to work. On the contrary, DCCA makes this 
problem naturally solvable (see \cite{VCSP, Vint, KV, V, VWsym}). Moreover, computations of basic
quantities related with nonlinear PDEs (symmetries, conservation laws, hamiltonian structures, 
B\"{a}cklund transformations, etc,etc) are essentially based on DCCA (see \cite{Vper, SKL, KK}). 
Secondary calculus, a natural language for the modern geometrical theory of nonlinear PDEs, is
the specialization of DCCA to diffieties. This explains the key importance of the DCCA-based methods, 
especially, cohomological ones for the theory of nonlinear PDEs. For further details about see
surveys \cite{Vfrom} and \cite{VWsym}.\newline

- {\it Graded differential geometry}.
As we have already pointed out the definition of  graded analogues of all objects of the standard 
differential geometry in terms of DCCA is identical to the non graded ones assuming that the latter 
are defined \emph{conceptually}. However, the search for conceptual definition of concrete 
quantities, the \emph{conceptualization problem}, could be a nontrivial task as  one can see from 
the previous discussion of tensors and integral forms. On the other hand, the recent history of 
formation of super-geometry when this opportunity was not taken into consideration illustrates the complications that could be otherwise avoided.
The necessity to systematically develop fundamentals of graded differential geometry on the basis
of DCCA is now urged by perspectives of important applications to both geometry and physics. 
The advantages that can be gained by putting problems in ordinary differential geometry into 
the graded framework are illustrated by the above ``conceptualization" of tensors in terms of iterated differential forms. Also, it is rather plausible that the problem of discretization of differential geometry 
for computer implementations is one of those that are naturally inscribed in this context. \newline

- {\it Bohr correspondence principle and observability in quantum physics}. 
Natural relations of quantum physics with the classical one that are partially expressed by the famous 
Bohr correspondence principle must be duly reflected in the formalising them mathematics. The 
commonly adopted von Neumann's proposal to appoint self-adjoint operators in Hilbert spaces
as observables in quantum mechanics manifestly violates the so generalised Bohr's principle.
Indeed, the very rich language of differential calculus that includes all the necessary for classical 
physics differential geometry has nothing in common with the poor and even rude language of 
Hilbert spaces. For instance, Hilbert spaces of functions defined on domains of different shapes 
and dimensions are isomorphic. Moreover, this language is not localisable in the space-time and 
also fails to work in QFT.

So, the problem of adequate mathematical formalisation of the observability mechanism in quantum
mechanics is still open as well as in QFT. The principal difficulty of this problem is that the physical 
factors ``responsible" for observability should be dully formalised and incorporated into this mechanism. 
The arguments of a smooth and natural correspondence with the language of classical physics suggest
to look for the details of this mechanism in graded differential geometry. The reader will find some more 
concrete ideas about in \cite{Vfrom} and \cite{conf}.  

Of a special interest in this list of problems is the long-standing problem of mathematically rigorous 
theory of integration by paths, which is crucial for modern QFT. It is plausible that the solution would 
be an analogue in secondary calculus of the cohomological theory discussed in subsection\,\ref{int}. 
\newline

\subsection{The language barrier}
One of the goals of this brief survey was to show with sufficiently simple examples that the language 
and methods of DCCA are both simplifying and  unifying. Not less important is that it reveals natural 
relations, which not infrequently remain hidden in the standard descriptive approach 
and, as a consequence, allows to better foresee the lines of the future developments by avoiding 
eventual misleading ideas and intuition that may come from the ``descriptiveness". And, finally, DCCA
makes possible a large expansion of methods of traditional ``differential mathematics" to new emerging
areas of mathematics and theoretical physics and also to not yet ``differentialised" domains of 
mathematics itself.

This is why we think that a systematic introduction of the language and methods of DCCA into the above indicated areas could have a notable positive effect.
Unfortunately, the language barriers slow down realisation in full of these, to our opinion, promising 
and intriguing possibilities. Another barriers are due to elevated diversity, complexity and dimension   
of the arising here problems, which require organisation of large scale research programs as it is 
common in experimental physics. In this connection a systematic introduction of fundamentals of 
DCCA into university courses could turn the situation to the best.\newline

{\bf Acknowledgements.} This paper presents an extended version of author's talk in the   
conference "Geometry of Jets and Fileds", 10-16 May, 2015, Bedlewo, Poland, dedicated to 60th
birthday of Janusz Grabowski.  The author takes this opportunity to express his warmest all the best
wishes to Janusz. He is very grateful to the organisers for the invitation and also to  M.M. Vinogradov 
and to L. Vitagliano who clarified him some points during the preparation of this paper.

\end{document}